\documentclass[12pt]{amsart}
\usepackage{amsmath}
\usepackage{amsfonts}
\usepackage{amssymb}
\usepackage{mathtools}
\usepackage{bbding}
\usepackage{stmaryrd}
\usepackage{txfonts}
\usepackage{booktabs}
\usepackage{graphicx}
\usepackage{graphics}
\usepackage{epsfig}
\usepackage{xypic}
\usepackage{tikz}
\usepackage{longtable}
\usepackage[all]{xy}
\usepackage{pgflibraryarrows}
\usepackage{pgflibrarysnakes}
\usepackage[shortlabels]{enumitem}

\usepackage{subfigure}
\usepackage{float}
    \usepackage[colorlinks,final,backref=page,hyperindex]{hyperref}

\usepackage{graphicx}
\usepackage{epstopdf}
\usepackage{epsfig}
\usepackage{pdfsync}    

\usepackage{xspace}

\newtheorem{theorem}{Theorem}[section]
\newtheorem{lemma}[theorem]{Lemma}
\newtheorem{corollary}[theorem]{Corollary}
\newtheorem{proposition}[theorem]{Proposition}
\theoremstyle{definition}
\newtheorem{definition}[theorem]{Definition}
\newtheorem{remark}[theorem]{Remark}
\newtheorem{example}{Example}[section]

\newcommand{\Der}{\operatorname{Der}}
\newcommand{\PDer}{\operatorname{PDer}}
\newcommand{\IDer}{\operatorname{IDer}}
\newcommand{\F}{\operatorname{Fix}}
\newcommand{\Id}{\textnormal{Id}_{A}}
\newcommand{\0}{\mathbf{0}_{A}}
\newcommand{\BA}{\mathbf{B}(A)}

\date{}

\usepackage{listings}
\usepackage{setspace}
\usepackage{xcolor}

\usepackage{listings}
\begin{document}
	
	\title[$(\odot,\vee)$-DERIVATIONS ON MV-ALGEBRAS]{$(\odot,\vee)$-DERIVATIONS ON MV-ALGEBRAS}
	
%
	
	\author{Xueting Zhao}
	\address{School of Mathematical Sciences, Shahe Campus, Beihang University, Beijing 102206, China}
	\email{xtzhao@buaa.edu.cn}
	
	\author{Aiping Gan}
	\address{School of Mathematics and Statistics,
		Jiangxi Normal University, Nanchang, Jiangxi 330022, P.R. China}
	\email{ganaiping78@163.com}
	
	\author{Yichuan Yang${}^\ast$}
	\address{School of Mathematical Sciences, Shahe Campus, Beihang University, Beijing 102206, China}
	\email{ycyang@buaa.edu.cn}
	\begin{abstract}
		Let $A$ be an MV-algebra. An $(\odot,\vee)$-derivation on $A$ is a map $d: A\to A$ satisfying:
		$d(x \odot y) = (d(x) \odot y) \vee(x \odot d(y))$
		for ~all $x, y \in A$. This paper initiates the study of $(\odot,\vee)$-derivations on MV-algebras. Several families of $(\odot,\vee)$-derivations on an MV-algebra are explicitly constructed to give realizations of the underlying lattice of an MV-algebra as lattices of $(\odot,\vee)$-derivations. Furthermore, $(\odot,\vee)$-derivations on a finite MV-chain are enumerated and the underlying lattice is described.

		
	\end{abstract}
	
	%
	%
	
		\maketitle
	\begin{center}
		\tableofcontents
	\end{center}
	
	\section{Introduction}\label{s01}
	
	
	~~~~The notion of derivation from analysis has been defined for various algebraic structures by extracting the Leibniz rule
	\begin{center}
		$\dfrac{d}{dx}(fg) = (\dfrac{d}{dx}(f))g + f(\dfrac{d}{dx}(g))$.
	\end{center}
	
	Derivations play an important role on describing the characteristics of prime rings \cite{prime}, and the multiplicative or additive commutativity of near rings \cite{near}, etc. A derivation in a prime ring $ (R, +, \cdot) $ is a map $ d $ : $  R\rightarrow R $ satisfying that for any $ x,y\in R $:
	$$ (1)~  d(x+y)=d(x)+d(y), \quad  (2)~ d(x\cdot y)=d(x)\cdot y+x\cdot d(y). $$
	The derivation on a lattice $ (L,\vee,\wedge) $ was defined by Szász \cite{lattice1975}, and was deeply investigated in \cite{lattice2001}, which is a map $ d $ : $ L\rightarrow L $ satisfying that for all $ x, y \in L $:
	$$  (i)~ d(x \vee y) = d(x) \vee d(y), \quad (ii)~ d(x \wedge y) = (d(x) \wedge y) \vee (x \wedge d(y)). $$
	The notion of derivations satisfying condition $ (ii) $ only was investigated by Xin and his coauthors \cite{lattice2008,lattice2012} with motivation from information science.
	In recent years the derivations have been defined and studied for BCI-algbras \cite{BCI}, BCC-algebras \cite{BCC2012,BCC2009}, BE-algebras \cite{BE}, and basic algebras \cite{BASIC}. Furthermore, the derivations on operator algebras were investigated by Bre\v{s}ar etc.\cite{radical,socle,local} which promoted the mathematical quantum mechanics and quantum field theory.
	
	An algebraic structure with a derivation is broadly called a differential algebra \cite{da}. In fact, differential algebra has found important applications in arithmetic geometry, logic and computational algebra, especially in the profound work of Wu on mechanical proof of geometric theorems \cite{DA1,DA2}. There are many instances of differential algebras, such as for fields \cite{field}, commutative algebras \cite{CA}, noncommutative algebras \cite{noncomm}, lattices \cite{DL}, and MV-algbras \cite{DMV2019}.
	
	The concept of derivations on MV-algebras was introduced by Alshehri \cite{DMV2010}: given an MV-algebra $ (M,\oplus,^{*} , 0) $, a derivation on $M$ is an operator (i.e, a map) $ d : M \rightarrow M $ such that $ d(x\odot y) = (d(x)\odot y) \oplus (x\odot d(y)), $ for all $ x, y \in M $, where $ x\odot y=(x^{*}\oplus y^{*})^{*} $. Furthermore, the different kinds of derivations on MV-algebras have been deeply investigated. Yazarli \cite{DMV2013} introduced the notions of symmetric bi-derivation, generalized derivation on MV-algebras. Then Wang, Davvaz and He \cite{DMV2017} studied additive derivations and their adjoint derivations to give a representation of MV-algebras. Recently, $\tau $-additive derivations on MV-algebras have been extended by Lu and Yang \cite{DMV2021}. Following these developments, we define the notion of $(\odot,\vee)$-derivations on $A$ satisfying $$d(x \odot y) = (d(x) \odot y) \vee(x \odot d(y))$$ for any $ x,y\in A $, where $ x\vee y=(x\odot y^{*})\oplus y $. 
     Our choice do not impose the extra ``union-preserving'' condition: $  d(x \vee y) = d(x)\vee d(y) $ and leads to several properties in this paper. Indeed as similar as \cite[Proposition 2.5]{DL}, a $(\odot,\vee)$-derivation with the ``union-preserving'' must be isotone.
	
	
	This paper initiates the study of $(\odot,\vee)$-derivations on MV-algebras. In Section \ref{s02}, we recall some necessary properties and examples of MV-algebras. In Section \ref{s03}, we introduce and study $(\odot,\vee)$-derivations on MV-algebras. After exploring a sufficient and necessary condition for an operator on a $n$-element MV-chain $L_{n}$ to be an $(\odot,\vee)$-derivation (Theorem \ref{thm03}), we show that the cardinality of the set of all $(\odot,\vee)$-derivations on $L_{n}$
	is  exactly $ \frac{(n-1)(n+2)}{2} $ (Theorem \ref{thm01}). 
	In Section \ref{s04}, the direct product of $(\odot,\vee)$-derivations is introduced. Let $\Omega$ be an index set, $\left\{A_{i}\right\}_{i \in \Omega}$ be a family of MV-algebras and $d_{i}$ be an operator of $A_{i}$ for each $i\in \Omega$, we prove that the direct product $\prod_{i \in \Omega} d_{i}$ of $ d_{i} $’s is an $(\odot,\vee)$-derivation (resp. a principal $(\odot,\vee)$-derivation) on $\prod_{i \in \Omega} A_{i}$
	if and only if $ d_{i} $ is an $(\odot,\vee)$-derivation (resp. a principal $(\odot,\vee)$-derivation) on $ A_{i} $ for each $ i \in \Omega$ (Theorem \ref{thm02}).
	In Section~\ref{s05}, we show that the set of $(\odot,\vee)$-derivations on a finite  MV-algebra has a natural lattice structure (Proposition \ref{pro12}) and we consider several lattice structure of  $(\odot,\vee)$-derivations which are isomorphic to the underlying lattice  $ \mathbf{L}(A) $ of an MV-algebra $A$ (Propositions \ref{pro14} and \ref{pro5.12}). We also describe the lattice structure of  $(\odot,\vee)$-derivations on finite MV-chains
	(Theorem \ref{thm04}).
	
	\textbf{Notations}. Throughout this paper, let $ |A| $ denote the cardinality of a set $ A $ and $ \mathbb{N}_{+} $ denote the set of all positive integers.
	
	\section{Preliminaries}\label{s02}
	
	~~~~In this section, we recall some necessary definitions and results about MV-algebras.
	
	\begin{definition}
		\cite[Definition 1.1.1]{MV4} An algebra $(A, \oplus, *, 0)$ of type $(2, 1, 0)$ is called an \textbf{MV-algebra} if it satisfies the following equations:
		
		(MV1) $ x\oplus(y\oplus z)=(x\oplus y)\oplus z;$
		
		(MV2) $ x\oplus y=y\oplus x;$
		
		(MV3) $ x\oplus 0=x;$
		
		(MV4) $  x^{**}=x;$
		
		(MV5) $ x\oplus  0^{*}= 0^{*};$
		
		(MV6) $  ( x^{*}\oplus y)^{*}\oplus y=( y^{*}\oplus x)^{*}\oplus x.$

	\end{definition}
	
	As usual, we shall denote an MV-algebra by its underlying carrier set.
	Note that all axioms of MV-algebras are equations, it follows by Birkhoff Theorem \cite[Theorem 11.9]{bu} that the class of all MV-algebras forms a variety.
	So the notions of isomorphism, subalgebra, congruence and direct product
	are just the particular cases of  the corresponding universal algebraic notions.
	
	\begin{example}\label{exa01}
		\textnormal{\cite{MV4}} 
		Let $ L =[0,1]$ be the real unit interval. Define 
		\begin{center}
			$ x \oplus y =
			\min\{1, x + y\} $~~and~~$ x^{*} = 1 - x $  for any $ x, y \in L$.	
		\end{center} 
		Then $ (L, \oplus, *, 0)$ is an MV-algebra.
		
		Let $ Q = [0, 1]\cap \mathbb{Q}$ and for each positive integer $ n \geq 2 $, let
		$$ L_{n} = \{0, \dfrac{1}{n-1}, \dfrac{2}{n-1}, \cdots, \dfrac{n-2}{n-1}, 1\}.$$ Then $ Q $ and the $n$-element subset $ L_{n} $ are subalgebras of $ L$ .
	\end{example}
	
	\begin{example}{\cite{MV1}}\label{exa04}
		Define the following sets of formal symbols:
		$$ \mathcal{C}_{0} = \{0, c, 2c, 3c, \cdots\},\quad  \mathcal{C}_{1} = \{1, c^{*},(2c)^{*},(3c)^{*},\cdots\}, $$
		where $ (kc)^{*} = 1 - kc $, and $ (kc)^{**} =((kc)^{*})^{*} = kc $ for any $ k \in \mathbb{N_{+}} $.
		
		Let $ + $ (respectively, $-$) be the ordinary sum (respectively, subtraction) between integers. We define the following binary operation $ \oplus $ on $ \mathcal{C}= \mathcal{C}_{0} \cup \mathcal{C}_{1} $:
		\begin{itemize}
			\item $ nc\oplus mc=(n+m)c $
			\item $ (nc)^{*}\oplus (mc)^{*}=1 $
			\item $ nc\oplus (mc)^{*}=(mc)^{*}\oplus nc=\left\{
			\begin{array}{ccc}
			1 & & {m\leq n}\\
			((m-n)c)^{*} & & {m>n}\\
			\end{array} 
			\right. $
		\end{itemize}
		
		Then $ (\mathcal{C}, \oplus, *,0) $ is an infinite MV-chain, and $0<c< 2c < 3c < \cdots < (n - 1)c < nc < \cdots < (nc)^{*} < ((n- 1)c)^{*} < \cdots < (3c)^{*} < (2c)^{*} < c^{*} < 1 $. MV-chains $Q $ and $ \mathcal{C} $ are not isomorphic, though they have the same countable cardinality.
		%
	\end{example}
	
	On every MV-algebra $A$, we define the constant $1$ and the operation
	$\odot$   as:
	$1 =0^{*}$ and $x\odot y=(x^{*}\oplus y^{*})^{*}.$    
	Then  for all $x, y\in A$, the following well-known
	properties hold \cite{MV4,MV5}:
	\begin{enumerate}
		\item [$\bullet$]    $(A, \odot, *, 1)$ is an MV-algebra;
		\item [$\bullet$]  $\ast$ is an
		isomorphism between $(A, \oplus, *, 0)$ and  $(A, \odot, *, 1)$;
		\item [$\bullet$]  $1^{*}=0, 1\oplus x=1$;
		\item [$\bullet$] $ x\oplus y=(x^{*}\odot y^{*})^{*}$;
		\item [$\bullet$] $x^{*} \oplus x = 1 $, 
		$ x\odot x^{*}=0 $.
	\end{enumerate}
	
	Let $ A $ be an MV-algebra. For any  $ x, y\in A $,
	define 
	$ x\leq y $ if and only if $x^{*}\oplus y=1$. Then
	$\leq$
	is a partial order on
	$A$, called \textbf{the natural order} of $A$  \cite{MV4}. Furthermore, the natural order determines a structure
	of bounded distributive lattice $ \mathbf{L}(A) $ on $A$, with $0$ and $1$
	are respectively the bottom and the top element, and
	$$x\vee y=(x\odot y^{*})\oplus y~~\textrm{and} ~~ x\wedge y=x\odot (x^{*}\oplus y).$$
	
	A linearly ordered MV-algebra  is called an\textbf{ MV-chain}. It is
	well-known that every $n$-element MV-chain is isormorphic to the MV-chain $L_{n}$ in Example \ref{exa01}.
	
	\begin{lemma}\label{lem01}
		\textnormal{\cite[Lemma 1.1.2]{MV4}}
		Let $ A $ be an MV-algebra and $ x, y \in A. $ Then the following statements are
		equivalent:
		\begin{enumerate}
			\item[$ (1) $] $x\leq y$;
			\item[$ (2) $] $x^{*}\oplus y=1$;
			\item[$ (3) $] $ x \odot y^{*} = 0 $;
			\item[$ (4) $] $ y = x \oplus (y \odot x^{*}) $;
			\item[$ (5) $] there is an element $ z \in A  $ such that $  x \oplus z = y $.
		\end{enumerate}
	\end{lemma}

	\begin{lemma}\label{lem02}	\textnormal{\cite{MV1,MV4}}
		Let $ A $ be an MV-algebra, and  $ x,y,z\in A $. Then the following statements hold:
		\begin{enumerate}
			\item[$(1)$] $ x\odot y\leq x\wedge y\leq x\leq x\vee y\leq x\oplus y $;
			
			\item[$(2)$] If $ x\oplus y=0 $, then $ x=y=0 $; If $ x\odot y=1 $, then $ x=y=1 $;
			
			\item[$(3)$] If $ x\leq y $, then $ x\vee z\leq y\vee z $, $ x\wedge z\leq y\wedge z $;
			
			\item[$(4)$] If $ x\leq y $, then $ x\oplus z\leq y\oplus z $, $ x\odot z\leq y\odot z $;
			
			\item[$(5)$] $ x\leq y $ if and only if $ y^{*}\leq x^{*} $;
			
			\item[$(6)$] $ x\odot (y\wedge z)=(x\odot y)\wedge (x\odot z) $;
			
			\item[$(7)$] $ x\odot (y\vee z)=(x\odot y)\vee (x\odot z) $;
			
			\item[$(8)$] $x \odot y\leq z$ if and only if $x\leq y^{*}\oplus z$.
		\end{enumerate}
	\end{lemma}
	
	\begin{lemma}\label{lem15}	\textnormal{\cite[Lemma 1.6.1]{MV4}}
		Let $ A $ be an MV-chain. For any $ x, y,z \in A$,
		\begin{enumerate}
			\item[$ (1) $] $ x\oplus y=x $ if and only if $ x=1 $ or $ y=0 $;
			\item[$ (2) $] If $x\odot y=x\odot z>0$, then $ y=z $.
		\end{enumerate} 
	\end{lemma}
	
	%
	%
	%
	
	\begin{example}\label{exa02}
		For any Boolean algebra  $(A, \vee, \wedge, -, 0, 1)$, the structure
		$(A, \vee, -, 0)$ is an MV-algebra, where $\vee, -$ and $0$ denote, respectively, the join,
		the complement and the smallest element in $A$.
	\end{example}
	
	Boolean algebras form a subvariety of the variety of MV-algebras. They
	are precisely the MV-algebras satisfying the additional equation $x \oplus x = x$. An element $a$ of $A$ is called \textbf{idempotent} if $a\oplus a=a$.
	Denote the set of all idempotent elements of $A$ by $\mathbf{B}(A)$, called \textbf{Boolean center of $ A $}. It is known that
	$ \textbf{B}(A) $ is a subalgebra of the MV-algebra $ A $, and a subalgebra $ B $ of $ A $ is a Boolean algebra if and only if
	$ B \subseteq \textbf{B}(A) $ \cite[Corollary 1.5.4]{MV4}.
	For convenience, we denote by $B_{n}$ the $n$-element Boolean algebra. It is clear that
	$B_{2}$ is exactly the  $2$-element MV-chain $L_{2}$.
	
	\begin{lemma}\label{lem03} \emph{\cite[Theorem 1.5.3]{MV4}} For every element $x$ in an MV-algebra $A$, the following conditions are equivalent:
		\begin{enumerate}
			\item[$(1)$]  $x\in \mathbf{B}(A)$;
			
			\item[$(2)$]  $x\oplus x=x$;
			
			\item[$(3)$]  $x\odot x=x$;
			
			\item[$(4)$]  $x^{*}\in \mathbf{B}(A)$;
			
			\item[$(5)$]  $x\oplus y=x\vee y$ for all $y \in A$;
			
			\item[$(6)$]  $x\odot y=x\wedge y$ for all $y \in A$.
		\end{enumerate}
	\end{lemma}
	
	
	
	\begin{definition}\label{def01}
		\textnormal{\cite{MV4}}	Let $ A $ be an MV-algebra and $ I $ be a  subset of $ A $. Then we say that $ I $ is an \textbf{ideal} if the following conditions are satisfied:
		\begin{enumerate}
			\item[$(1)$] $ 0\in I $;
			
			\item[$(2)$] $ x,y\in I $ imply $ x\oplus y\in I $;
			
			\item[$(3)$] $ x\in I $ and $ y\leq x $ imply $ y\in I $.	
		\end{enumerate}
	\end{definition}
	\begin{definition}\label{def04}
		\textnormal{\cite{MV4}}	Let $ A $ be a lattice and $ I $ be a subset of $ A $. Then we say that $ I $ is a \textbf{lattice ideal} if the following conditions are satisfied:
		\begin{enumerate}
			\item[$(1)$] $ 0\in I $;
			
			\item[$(2)$] $ x,y\in I $ imply $ x\vee y\in I $;
			
			\item[$(3)$] $ x\in I $ and $ y\leq x $ imply $ y\in I $.	
		\end{enumerate}
	\end{definition}
	That is, a lattice ideal of an MV-algebra $ A $ is the notion of ideal in the underlying lattice $ (A, \wedge, \vee) $\cite[Proposition 1.1.5]{MV4}. It can easily be verified that an ideal is a lattice ideal but the opposition is not necessarily the case. The next lemma gives the representation of a finite MV-algebra: 
	
	\begin{lemma}\label{lem12} \textnormal{\cite[Propostion 3.6.5]{MV4}}
		An MV-algebra $A$ is finite if and only if $A$ is isomorphic to a finite product of finite chains, in symbols,
		$$A \cong L_{d_{1}} \times \cdots \times L_{d_{u}},$$ 
		for some integers $2 \leq d_{1} \leq d_{2} \leq \ldots \leq d_{u}$.
		This representation is unique, up to the ordering of factors.
	\end{lemma}
	
	Finally, we list the famous Chang's Subdirect Representation Theorem, stating that if an equation holds in all totally ordered MV-algebras, then the equation holds in all MV-algebras.
	
	\begin{lemma}\label{lem14} \textnormal{\cite[Theorem 1.3.3]{MV4}}
		Every nontrivial MV-algebra is a subdirect product of MV-chains.
	\end{lemma}
	
	%
	
	\section{$(\odot,\vee)$-derivations on MV-algebras}\label{s03}
	
	~~~~In this section, we introduce $(\odot,\vee)$-derivations on MV-algebras, and characterize some properties about $(\odot,\vee)$-derivations, such as
	isotonicity and idempotency.
	Also, we enumerate
	the cardinality of $(\odot,\vee)$-derivations on finite MV-chains. 
	
	\subsection{Basic properties of $(\odot,\vee)$-derivations on MV-algebras}\label{sec3.1}
	
	\begin{definition}
		\textnormal{Let $ A $ be an MV-algebra. A map $d : A\rightarrow A $  is called an \textbf{$(\odot,\vee)$-derivation on $A$} if it satisfies the  equation:
			\begin{equation}\label{equ01}
			d(x \odot y) = (d(x) \odot y) \vee(x \odot d(y)) \quad
			for ~all ~x, y \in A.	
			\end{equation}
		}
	\end{definition}
	
	It is easy to check that the identity map $\textnormal{Id}_{A}$ and the zero map $\textbf{0}_{A}$ are simple examples of $(\odot,\vee)$-derivations on an MV-algebra $A$, where
	\begin{center}
		$\textnormal{Id}_{A}(x) = x$ \quad  and \quad $\textbf{0}_{A}(x)=0$ \quad
		for  any  $x \in A$.	
	\end{center}
	
	Also, for a given $a\in A$, define the map $d_{a}: A\to A$ by
	\begin{center}
		$d_{a}(x):= a\odot x$ \quad for all $x\in A.$
	\end{center}
	
	Then $d_{a}$ is an $(\odot,\vee)$-derivation, called a \textbf{principal $(\odot,\vee)$-derivation}.
	Both  $\textnormal{Id}_{A}$ and  $\textbf{0}_{A}$ are principlal $(\odot,\vee)$-derivations, since
	$\textnormal{Id}_{A}=d_{1}$ and $\textbf{0}_{A}=d_{0}$.

	Denote the set of all $(\odot,\vee)$-derivations on $ A$ by $\Der(A)$;  and the set of
	all the  principal $(\odot,\vee)$-derivations on $A$  by $ \PDer(A) $, that is
	$ \PDer(A)=\{d_{a} ~|~a\in A\}. $

	\begin{remark}\label{rem:3.2}
		\begin{enumerate}
			\item[$(1)$] It is clear that Eq.$\eqref{equ01}$ holds when $x=y=1$, where $ d(1)=d(1)\odot 1 $.
			
			\item[$(2)$] Adapting the classical terminology of differential algebras, we also call a derivation a differential operator. More generally, we also call a map $f:A\to A$ an operator even though there is no linearity involved.
			
			\item[$(3)$] Note that in \cite{DMV2010,DMV2017}, an $(\odot,\oplus)$-derivation on an MV-algebra $A$ is defined to be a map satisfying $ d(x \odot y) = (d(x) \odot y) \oplus(x \odot d(y))$ for all $x, y\in A$. In this paper, we use ``$ \vee $" instead of ``$ \oplus $". Our choice of this notation has its motivation from certain asymmetry of ``$ \vee $" and ``$ \odot $", and already leads to some properties as displayed in Proposition \ref{pro01}. 

			\item[$(4)$]  It is natural to consider a $(\oplus,\wedge)$-derivation which is dual to the $(\odot,\vee)$-derivation on an MV-algebra $A$: $d(x \oplus y)=(d(x)\oplus y)\wedge(x\oplus d(y))$ for all $ x, y\in A$. If this condition is taken, then the study should be completely parallel to the study of Eq. \eqref{equ01} due to the symmetry of the operations ``$ \vee $'' and ``$ \wedge $'', ``$ \odot $'' and ``$ \oplus $'' in the definition of an MV-algebra. Furthermore, if a map $d$ is both an $(\odot,\vee)$-derivation and 
			a $(\oplus,\wedge)$-derivation,
			then $ d=\textnormal{Id}_{A} $ (see Proposition \ref{pro03}). 
		\end{enumerate}
	\end{remark}

	\begin{proposition}\label{pro01}
		Let $A$ be an MV-algebra, $x, y \in A$ and $d\in \Der(A)$. Then  for any positive integer $n$, the following statements hold:
		\begin{enumerate}
			\item[$(1)$]  $d(0) = 0$.
			
			\item[$(2)$]	 $d(x^{n}) = x^{n-1}\odot d(x)$, where $x^{0}=1$, $x^{n}=\overbrace{x\odot x \odot \cdots \odot x}^{n}$.
			
			\item[$(3)$] $d(x)\odot x^{*}= x \odot d(x^{*})  = 0$.
			
			\item[$(4)$]	 $d(x) \leq x$.
			
			\item[$(5)$]	 $d(x) = d(x) \vee(x \odot d(1))$ and so $x \odot d(1)\leq d(x)$.
			
			\item[$(6)$]	 $d(x^{*}) \leq x^{*}\leq (d(x))^{*}$.
			
			\item[$(7)$] $d(x)\odot d(y) \leq d(x \odot y) \leq d(x) \vee d(y) \leq d(x) \oplus d(y)$.
			
			\item[$(8)$]	 $(d(x))^{n}\leq d(x^{n})$.
			
			\item[$(9)$] If $I$ is a downset of $A$, then $d(I)\subseteq I$, where $d(I)=\{d(x) | x\in I\}$.
			
			\item[$(10)$]
			If $ y\leq x $ and $ d(x)=x $, then $ d(y)=y $.
		\end{enumerate}
	\end{proposition}
	
	\begin{proof}
		$(1)$ Putting $x=y=0$ in Eq.\eqref{equ01}, we immediately have $d(0) = d(0\odot0) =
		(d(0)\odot 0)\vee(0 \odot d(0)) =0 $.
		\vspace{0.20cm}
		
		\noindent$(2)$ We prove $d(x^{n}) = x^{n-1}\odot d(x)$ by induction on $n$. First,
		it is clear that $d(x^{1}) =d(x)=1 \odot d(x)= x^{1-1}\odot d(x)$. For $n=2,$
		putting $x = y$ in Eq.\eqref{equ01}, we get $d(x^{2}) = d(x\odot x) = (d(x) \odot x)\vee 
		(x\odot d(x)) = x \odot d(x)$. 
		
		Now assume that $d(x^{n}) = x^{n-1}\odot d(x)$. By Eq.\eqref{equ01}, we have $d(x^{n+1})=d(x^{n}\odot x) =(d(x^{n})\odot x)\vee (x^{n}\odot d(x)) =(x^{n-1}\odot d(x)\odot x)\vee (x^{n}\odot d(x))=x^{n}\odot d(x)$, and so $(2)$ holds.
		\vspace{0.20cm}
		
		\noindent$(3)$ Since $x\odot x^{*}=0$, by Item $(1)$ it follows that $0=d(0)=d(x\odot x^{*})=(d(x)\odot x^{*})\vee (x\odot d(x^{*}))$. So $d(x)\odot x^{*}=0$ and $ x \odot d(x^{*})=0$. 
		\vspace{0.20cm}
		
		\noindent$(4)$ Since  $d(x)\odot x^{*}=0$ by Item $(3)$, it follows immediately by Lemma \ref{lem01} that $d(x)\leq x$.
		\vspace{0.20cm}
		
		\noindent$(5)$ By Eq.\eqref{equ01} we have $d(x) = d(x\odot 1) = (d(x) \odot 1)\vee(x\odot d(1)) =
		d(x) \vee(x \odot d(1))$. So $x \odot d(1) \leq d(x)$.
		\vspace{0.20cm}
		
		\noindent$(6)$ We have $d(x^{*}) \leq x^{*}$ and $ d(x)\leq x $   by Item $(4)$. Thus  $ x^{*}\leq (d(x))^{*} $ by  Lemma \ref{lem02} $(5)$.
		\vspace{0.20cm}
		
		\noindent$(7)$ By Item $(4)$ and Lemma \ref{lem02} $(4)$, we have $ d(x)\odot d(y) \leq x\odot d(y)$ and $d(x)\odot d(y) \leq d(x)\odot y $. So $d(x)\odot d(y) \leq (d(x)\odot y)\vee(x\odot d(y)) = d(x\odot y)$.
		Furthermore, it follows from Lemma \ref{lem02} $(1)$ that $ d(x) \odot y \leq d(x),x \odot d(y) \leq d(y)$. So $d(x\odot y)=(d(x)\odot y) \vee(x \odot d(y))   \leq d(x) \vee d(y)$. Finally, we  get $d(x)\vee d(y)\leq d(x) \oplus d(y)$  by Lemma \ref{lem02} $(1)$.
		\vspace{0.20cm}
		
		\noindent$(8)$ By Item (2), we have $ d(x^{n}) = x^{n-1}\odot d(x) $. Since $ d(x)\leq x $, it follows by Lemma \ref{lem02} (4) that $ (d(x))^{n-1}\leq x^{n-1} $ and then $(d(x))^{n}=(d(x))^{n-1}\odot d(x)\leq x^{n-1}\odot d(x)=d(x^{n})$. 
		\vspace{0.20cm}
		
		\noindent$(9)$ Let $I$ be a downset of $A$ and $y\in d(I)$. Then there exists $a\in I$ such that $y=d(a)$. Since $d(a)\leq a$  by Item $(4)$, we have $y=d(a)\in I$ by Definition \ref{def01}. Thus $d(I)\subseteq I$.
		\vspace{0.20cm}
		
		\noindent$(10)$	If $ y\leq x $ and $d(x)=x$, then
		\begin{eqnarray*}
			&&d(y)=d(x\wedge y)=d(x\odot(x^{*}\oplus y))\\
			&=&(d(x)\odot(x^{*}\oplus y))\vee (x\odot d(x^{*}\oplus y))\\
			&=&(x\odot(x^{*}\oplus y))\vee (x\odot d(x^{*}\oplus y))\\
			&=&x\odot(x^{*}\oplus y)\\
			&=&x\wedge y\\
			&=&y,
		\end{eqnarray*}	
		and so we get  $ d(y)=y$.
	\end{proof}

	It is known that if $d$
	is a derivation on a lattice $L$, then $ d=\textnormal{Id}_{L}$ iff $d$ is injective iff $d$ is surjective 
	\cite[Theorem 3.17]{lattice2008}. In Proposition \ref{pro03}, we will show that 
	if $d$
	is an $(\odot,\vee)$-derivation on an MV-algebra $A$, then $ d=\textnormal{Id}_{A}$ iff $d$ is surjective. However, $d$ is injective may not imply that $ d=\textnormal{Id}_{A}$ (see Remark \ref{remark3.5}).
	
	\begin{proposition}\label{pro03}
		Let $A$ be an MV-algebra and $d\in \Der(A)$.	 Then the following statements are equivalent:
		\begin{enumerate}
			\item[$(1)$] $ d=\textnormal{Id}_{A} $;
			\item[$(2)$] $ d(1)=1 $;
			\item[$(3)$] $ d(a)=1 $ for some $ a\in A $;
			\item[$(4)$] $ d $ is surjective;
			\item[$(5)$] $d$ is a $(\oplus,\wedge)$-derivation, i.e., $d$ satisfies the condition:
			$d(x \oplus y)=(d(x)\oplus y)\wedge(x\oplus d(y))$ for all $ x, y\in A$.
		\end{enumerate}
		
	\end{proposition}
	\begin{proof}
		It is clear that $(1)\Rightarrow (2)\Rightarrow (3)$, and $(1)\Rightarrow (4)\Rightarrow (3)$ by the property of $ \textnormal{Id}_{A} $.
		\vspace{0.20cm}
		
		\noindent$(2)\Rightarrow (1)$. Assume that $d(1)=1$. Then by Proposition \ref{pro01} $(10)$ we have that $ d(x)=x $ for all $ x\in A $. Thus $ d=\textnormal{Id}_{A}$,
		Item $(1)$ holds.
		\vspace{0.20cm}
		
		\noindent$(3)\Rightarrow (2)$. Assume that $d(a)=1$ for some $a\in A$. By Proposition \ref{pro01} $ (4) $, we have $ 1=d(a)\leq a $, and so $a=1 $. Thus $ d(1)=1 $, Item $(2)$ holds.
		\vspace{0.20cm}
		
		\noindent$(1)\Rightarrow (5)$. Assume that $ d=\textnormal{Id}_{A} $. Then $ d(x\oplus y)=x\oplus y=(x\oplus y)\wedge(x\oplus y)=(d(x)\oplus y)\wedge(x\oplus d(y)) $ for all $x, y\in A$, and thus Item $(5)$ holds.
		\vspace{0.20cm}
		
		\noindent$(5)\Rightarrow (2)$. Assume that $ d $ is a $(\oplus,\wedge)$-derivation. Then $ d(1)=d(1\oplus1)=(d(1)\oplus 1)\wedge(1\oplus d(1))=1 $, and so Item $(2)$ holds.
	\end{proof}
	
	\begin{remark} \label{remark3.5}	
		Let $A$ be an MV-algebra and $d\in \Der(A)$. Generally,	  $d\neq\textnormal{Id}_{A}$ if 
		$d$ is  injective. 
		For example, let  $\mathcal{C} $ be the infinite MV-chain in Example \ref{exa04}. 
		Define an operator $ d $ on  $\mathcal{C} $ by
		$$d(x):=
		\begin{cases}
		x\odot c^{*},  & \textrm{if}~ x\in \mathcal{C}_{1} \\
		x,  & \textrm{if}~ x\in \mathcal{C}_{0}
		\end{cases}
		$$
		
		\textbf{Claim} $(1)$:  $d\in \Der(\mathcal{C}) $.
		Indeed, let $x,y\in\mathcal{C}$. Consider the following cases:
		
		Case $(i)$:	$x, y\in \mathcal{C}_{1}$. Then $ d(x\odot y)=(x\odot y)\odot c^{*}=(d(x)\odot y)\vee(x\odot d(y))  $. 
		
		Case $(ii)$: $x, y\in \mathcal{C}_{0} $. Then
		$ d(x\odot y)=x\odot y=(d(x)\odot y)\vee(x\odot d(y))$. 
		
		Case $(iii)$: $x\in\mathcal{C}_{1}, y\in\mathcal{C}_{0}$, let $ x=(mc)^{*}, y=nc$, where $m, n\in \mathbb{N_{+}}$. Then $d(x)=(mc)^{*} \odot c^{*}=((m+1)c)^{*}$ and $d(y)=y$.

		If $ m\geq n $, then 
		$ d(x\odot y)=d((mc)^{*}\odot nc)=d(0)=0=(((m+1)c)^{*}\odot nc)\vee((mc)^{*}\odot nc)=(d(x)\odot y)\vee (x\odot d(y)) $. If $ m<n $,
		then 
		$$d(x)\odot y=((m+1)c)^{*}\odot nc=
		\begin{cases}
		0,  & \textrm{if}~ m+1=n \\
		(n-m-1)c,  & \textrm{if}~ m+1<n
		\end{cases}
		$$ 
		It follows that $d(x)\odot y< (n-m)c=x\odot d(y)$, and thus $d(x\odot y)=d((mc)^{*}\odot nc)=d((n-m)c)=(n-m)c=(d(x)\odot y)\vee(x\odot d(y))  $. 
		
		Case $(iv)$: $x\in\mathcal{C}_{0}, y\in\mathcal{C}_{1}$. Similarly,  we can obtain that  $ d(x\odot y)=(d(x)\odot y)\vee(x\odot d(y))$.
		
		Summarizing the above arguments, we get $d\in \Der(\mathcal{C}) $.
		
		\textbf{Claim} $(2)$: $ d $ is injective. Indeed, let $x,y\in\mathcal{C}$ and $x\neq y$. If $x, y\in \mathcal{C}_{1}$, say $ x=(mc)^{*},y=(nc)^{*} $, where $m, n$ are positive integers and $ m\neq n $, then $d(x)=(mc)^{*} \odot c^{*}=((m+1)c)^{*}\neq ((n+1)c)^{*}=(nc)^{*} \odot c^{*}=d(y)$. 	
		If $x, y\in \mathcal{C}_{0}$, then $ d(x)=x\neq y=d(y) $.
		
		If $x\in\mathcal{C}_{1}, y\in\mathcal{C}_{0}$ or $y\in\mathcal{C}_{1}, x\in\mathcal{C}_{0}$, say $x\in\mathcal{C}_{1}, y\in\mathcal{C}_{0}$,
		then by the definition of $d$, we have 
		$d(x)\in\mathcal{C}_{1},  d(y)\in\mathcal{C}_{0}$, so $ d(x)\neq d(y) $ since $\mathcal{C}_{0}\cap \mathcal{C}_{1}=\emptyset$. Thus $d$ is injective.
		
		However,  $d\neq \textnormal{Id}_{\mathcal{C}} $ since $d(1)=c^{*}\neq 1$.	
	\end{remark}
	
	Let $A$ be an MV-algebra and $d\in \Der(A)$.
	From  Remark \ref{remark3.5}, we see that $d(a)$ may not lie in $\text{\textbf{B}}(A)$  if $a\in \text{\textbf{B}}(A)$. 
	In what follows, some properties of $(\odot,\vee)$-derivations   related to Boolean center $ \text{\textbf{B}}(A) $ of an MV-algbera $A$ are given.

	\begin{proposition}\label{pro11}
		Let $A$ be an MV-algebra and $d\in \operatorname{Der}(A)$.	 Then for all $x, y \in \textnormal{\textbf{B}}(A)$,
		the following statements hold:
		\begin{enumerate}
			\item[$(1)$] $d(x \wedge y) = (d(x) \wedge y) \vee (x \wedge d(y))$.
			
			\item[$(2)$] $d(x) = x \odot d(x)$.
		\end{enumerate}
	\end{proposition}
	\begin{proof} 
		
		\noindent$(1)$ By Lemma \ref{lem03} $(6)$, we have 
		$d(x\wedge y) = d(x\odot y) = (d(x)\odot y)\vee (x\odot d(y)) = (d(x)\wedge y)\vee (x\wedge d(y))$.
		\vspace{0.20cm}
		
		\noindent$(2)$ Since  $ x\odot x=x $, we have $d(x)=d(x\odot x)=x\odot d(x) $ by Proposition \ref{pro01} $(2)$.
	\end{proof} 
	
	\begin{corollary}\label{cor04}
		If an MV-algebra $A$ is a Boolean algebra, then $d$ is an $(\odot,\vee)$-derivation on $A$ if and only if $d$ is a derivation on the lattice $(A,\vee,\wedge)$.
	\end{corollary}
	\begin{proof}
		It follows immediately by Proposition \ref{pro11} and Lemma \ref{lem03}.
	\end{proof}

	Note that $ d(d(a))$ may not equal $d(a) $  if $ a\in \textnormal{\textbf{B}}(A) $. For example, in Remark \ref{remark3.5}, we have $1\in \textnormal{\textbf{B}}(\mathcal{C})$ but $d(d(1))=d(c^{*})=c^{*}\odot c^{*}=(2c)^{*}\neq c^{*}=d(1)$.
	Proposition \ref{pro3.7} tells us  that $d(d(a)) = d(a)$ if $d(a) \in \textnormal{\textbf{B}}(A)$.

	
	\begin{proposition}\label{pro3.7}
		Let $A$ be an MV-algebra, $d\in \Der(A)$	and $a\in A$. If $d(a) \in \textnormal{\textbf{B}}(A)$,
		then $d(d(a)) = d(a)$.
	\end{proposition}
	
	\begin{proof}
		Assume that $d\in \Der(A)$, $a\in A$	with $d(a) \in \textbf{B}(A)$. Then  $d(a)=d(a)\odot d(a) \leq d(a\odot a)=a\odot d(a)\leq d(a)$ by Proposition \ref{pro01} $(8)$ and
		Lemma \ref{lem02} $(1)$. Thus $ d(a)=a\odot d(a) $, and therefore $ d(d(a))=d(a\odot d(a))=(d(a)\odot d(a))\vee (a\odot d(d(a)))=d(a)\vee (a\odot d(d(a))) $ by Eq.~\eqref{equ01}. Consequently, we get 
		$d(a) \leq d(d(a))$. Also, we have $d(d(a)) \leq d(a)$ by Proposition \ref{pro01} (4).
		Hence $d(d(a)) = d(a)$. 
	\end{proof}

\subsection{$(\odot,\vee)$-derivations on MV-chains}\label{sec3.2}
	\ 
	\newline 
	\indent
	
	~~~~In this subsection we will determine the cardinality of $\Der (A)$ when $A$ is a finite MV-chain.
	Let $n\geq 2$ be a positive integer. Recall that 
	every $n$-element MV-chain is isomorphic to the MV-chain $L_{n}$, where  $L_{n}$ is given in Example \ref{exa01}. 
	
	\begin{remark}\label{rem:3.9}
		In $L_{n}$, $\frac{n-m-1}{n-1}=(\frac{n-2}{n-1})^{m}$ for each $m\in \{1, 2, \cdots, n-1\}$. That is to say,
		for any $x\in L_{n}\backslash \{1\}$,  $ x $ can be expressed as a power of $ \frac{n-2}{n-1} $.	
	\end{remark}

	\begin{theorem}\label{thm03}
		Let  $d$ be an operator on  $L_{n}$ and $ v=\frac{n-2}{n-1} $. Suppose that $ d(v)\leq v$. Then $d\in \Der(L_{n})$ if and only if $d$ satisfies the following conditions:
		\begin{enumerate}
			\item[$ (1) $]  $d(v^{m})=v^{m-1}\odot d(v)$ for each $m\in \{1, 2, \cdots, n-1\}$;
			\item[$ (2) $]  $v\odot d(1)\leq d(v) $.
		\end{enumerate}
		
	\end{theorem}
	\begin{proof}
		
		If $d\in \Der(L_{n})$, then for each $m\in \{1, 2, \cdots, n-1\}$, we have $d(v^{m})=v^{m-1}\odot d(v)$ 
		by Proposition \ref{pro01} $(2)$, and $v\odot d(1)\leq d(v \odot 1)=d(v) $ 
		by Proposition \ref{pro01} $(5)$.
		Thus $d$ satisfies the conditions $ (1) $ and $ (2) $.
		
		Conversely, suppose that $d$ satisfies the conditions $ (1) $ and $ (2) $. Let $x, y\in L_{n}$. By Remark \ref{rem:3.2} $(1)$,
		we can assume that $x\neq 1$ or $y\neq 1$ and distinguish the following cases:
		
		If $x\neq 1$ and $y \neq 1$, then $ x=v^{k} $ and $ y=v^{l} $ for some $k, l\in \{1, 2, \cdots, n-1\}$. By the condition (1), we  get $ d(x\odot y)=d(v^{k}\odot v^{l})=v^{k+l-1}\odot d(v)=((v^{k-1}\odot d(v))\odot v^{l})\vee (v^{k}\odot (v^{l-1}\odot d(v)))=(d(x)\odot y)\vee(x\odot d(y))$. 
		
		If $x=1$ or $y=1$ (but not both), say $x \neq 1$ and $y = 1$, then 
		$ x=v^{k} $ for some $k\in \{1, 2, \cdots, n-1\}$.
		By the condition (1), we have $d(x)=d(x\odot1)= d(v^{k})=v^{k-1}\odot d(v) $. Also, we have $ x\odot d(1)=v^{k-1}\odot v\odot d(1)\leq v^{k-1}\odot d(v)=d(x)\odot 1$ by condition (2). Thus we have derived that $d(x\odot 1)=  d(x)=d(x)\odot 1=(d(x)\odot 1)\vee(x\odot d(1))$.
		
		Therefore, we conclude that $d\in \Der(L_{n})$.	
	\end{proof}

	From Theorem \ref{thm03}, we see that if  $d\in \Der(L_{n})$, then for any $x\in L_{n}$ with $x<\frac{n-2}{n-1} $, 
	$d(x)$ is determined by the value  $d(\frac{n-2}{n-1})$.
	However, if $L$ is
	an infinite MV-chain with an anti-atom $v ~($i.e, $v$ is the maximum element in $L\backslash \{1\} )$ and $d\in \Der(L)$, then for any $x<v $, 
	$d(x)$ may not be determined by the value  $d(v)$.
	For example, let $ \mathcal{C} $  be the MV-chain in Example \ref{exa04}.
	Then $c^{*}$ is the anti-atom of $ \mathcal{C} $.
	Define operators $ d $ and $ d' $ on $ \mathcal{C} $ as follows:
	$$d(x):=
	\begin{cases}
	x\odot c^{*},  & \textrm{if}~ x\in \mathcal{C}_{1} \\
	x,  & \textrm{if}~ x\in \mathcal{C}_{0}
	\end{cases}\quad and \quad
	d'(x):=x\odot c^{*}
	$$
	Then  $d\in \Der(\mathcal{C})$ by Remark \ref{remark3.5} and $d'$ is a principal $(\odot, \vee)$-derivation.
	Furthermore, $d(c^{*})=d'(c^{*})$
	but $ d\neq d' $ since $d(c)=c\neq 0=d'(c)$.

	\begin{theorem}\label{thm01}
		Let $n\geq 2 $ be a positive integer. Then $|\Der(L_{n})|= \frac{(n-1)(n+2)}{2}$.
	\end{theorem}
	
	\begin{proof}
		Assume that $d\in \Der(L_{n})$ and denote $\frac{n-2}{n-1}$ by $v$. Then
		$d(v)\leq v$ by Proposition \ref{pro01} $(4)$, and so
		$ d(v)=\frac{i}{n-1} $ for some $i \in \{0, 1, 2, \cdots, n-2\} $.
		For  any $x\in L_{n}$ with 
		$x<v $, 
		$d(x)$ is determined by the value  $d(v)$ by Theorem \ref{thm03}.
		
		Now consider the value $d(1)$.
		By the condition $ (2) $ of Theorem \ref{thm03}, we have 
		\begin{equation}\label{equ2}
		v\odot d(1)=\frac{n-2}{n-1}\odot d(1)\leq d(v)= \frac{i}{n-1}. 
		\end{equation}
		Notice that for all $k, l\in \{0, 1, 2, \cdots, n-2\}$, we have $\frac{k}{n-1}\odot \frac{l}{n-1}=\max \{0, \frac{k+l}{n-1} -1 \}$.
		Eq. \eqref{equ2} implies that $ d(1)\leq \frac{i+1}{n-1} $. So $d(1)$ has $i+2$ choices.
		
		Summarizing the above arguments, we get
		\[|\Der(L_{n})|=\sum_{i=0}^{n-2}(i+2)=2+3+\cdots+n=\frac{(n-1)(n+2)}{2}.\qedhere\]	
	\end{proof}
	
	By Theorem \ref{thm01}, we obtain 
	$|\Der( L_{2} )|=\frac{(2-1)(2+2)}{2}=2$ and
	$|\Der( L_{3} )|=\frac{(3-1)(3+2)}{2}=5$. Thus
	$\Der (L_{2})=\{\textnormal{Id}_{L_{2}},
	\textbf{0}_{L_{2}} \}$. Let $A$ be an MV-algebra.
	In what follows, we will show that  $|\Der(A )|=2$ iff $A$ is isomorphic to $L_{2}$; and $|\Der(A )|=5$ iff $A$ is isomorphic to $L_{3}$. For this purpose, we first give a family of derivations on $A$.

	\begin{proposition}\label{pro08}
		Let $A$ be an MV-algebra and $d\in \operatorname{Der}(A)$. Let $u \in A$ be given with $u \leq d(1)$ and
		define an operator $ d^{u} $ on $ A $ by
		$$d^{u}(x):=
		\begin{cases}
		u& \textrm{if} ~~ x = 1\\
		d(x)& \textrm{otherwise}
		\end{cases}$$
		Then $ d^{u} $ is also in $\Der(A)$. 
	\end{proposition}
	\begin{proof}
		Let $x,y \in A$. By Remark \ref{rem:3.2} $(1)$,
		we can assume that $x\neq 1$ or $y\neq 1$.
		
		If $x\neq 1$ and $y \neq 1$, then
		$d^{u}(x)=d(x)$, $d^{u}(y)=d(y)$ and $x \odot y \in A\setminus\{1\}$, which implies that
		$d^{u}(x \odot y) = d(x \odot y) = (d(x) \odot y) \vee (x \odot d(y)) = (d^{u}(x) \odot y) \vee (x \odot d^{u}(y))$.
		
		If $x=1$ or $y=1$ (but not both), say
		$x \neq 1$ and $y = 1$, then since $d^{u}(1) = u \leq d(1)$, we have $x\odot d^{u}(1) \leq x \odot d(1) \leq d(x)$ by Proposition \ref{pro01} $(4)$ and so
		$$d^{u}(x\odot y) = d^{u}(x) = d(x) = d(x) \vee (x \odot d^{u}(1)) = (d^{u}(x) \odot y) \vee (x \odot d^{u}(y)).$$
		
		
		
		Thus we conclude that $ d^{u} $ is in $\Der(A)$.
	\end{proof}
	
	\begin{corollary}\label{pro09}
		Let $A$ be an MV-algebra, and $u\in A$. Define operators $\chi^{(u)}$ as follows:
		$$\chi^{(u)}(x):=
		\begin{cases}
		u,  & \text{if}~ x=1 \\
		x,  & \text{otherwise}.
		\end{cases}\quad  \quad$$
		Then $\chi^{(u)}\in \Der(A)$.
	\end{corollary}
	\begin{proof}
		Since $\textnormal{Id}_{A}\in \Der(A)$ and $u\leq 1=\textnormal{Id}_{A}(1)$, we have
		$\chi^{(u)}=(\textnormal{Id}_{A})^{u}\in \Der(A)$ by Proposition \ref{pro08}.
	\end{proof}
	
	\begin{lemma}\label{lem09}
		Let $A$ be an MV-algebra. Then the following statements hold:
		\begin{enumerate}
			\item[$ (1) $] $\chi^{(0)}\neq d$ for any $d\in \Der(A)$ with $d(1)\neq 0$. In particular, $\chi^{(0)}\neq \chi^{(u)}$ and $\chi^{(0)}\neq d_{u}$ for any $u\in A\backslash \{0\}$.
			
			\item[$ (2) $] If $|A|\geq 3$, then $\chi^{(u)}\neq d_{v}$ for any  $u, v\in A\backslash\{0, 1\}$.
		\end{enumerate}
	\end{lemma}
	\begin{proof}
		$ (1) $   Since $\chi^{(0)}(1)=0$, it follows that $\chi^{(0)}\neq d$ for any $d\in \Der(A)$ with $d(1)\neq 0$, which implies that
		$\chi^{(0)}\neq \chi^{(u)}$ and $\chi^{(0)}\neq d_{u}$ for any $u\in A\backslash \{0\}$, since  $\chi^{(u)}(1)=d_{u}(1)=u\neq 0$.
		\vspace{0.20cm}

		\noindent$ (2) $ Assume that $|A|\geq 3$ and let $u, v\in A\backslash\{0, 1\}$.
		Then $u^{*}, v^{*}\in A\backslash\{0, 1\}$. 
		
		If $ u\neq v $, then $ \chi^{(u)}\neq d_{v} $, since
		$\chi^{(u)}(1)=u  $  $\neq v= d_{v}(1)$. 
		If $ u=v $, then $ \chi^{(u)}\neq d_{u} $, since 
		$\chi^{(u)}(u^{*})=u^{*}\neq0=u\odot u^{*}=d_{u}(u^{*})$.
	\end{proof}
	
	\begin{corollary}\label{cor02}
		Let $A$ be an MV-algebra. Then the following statements hold:
		\begin{enumerate}
			\item[$ (1) $]  If $|A|\geq 3$, then $ |\Der(A)|\geq 5 $ .
			
			\item[$ (2) $] If $|A|\geq 4$, then $ |\Der(A)|\geq 7 $.
			
			\item[$ (3) $]  If $|A|\geq 5$, then $ |\Der(A)|\geq 13 $.
		\end{enumerate}		
	\end{corollary}
	\begin{proof}
		$ (1) $ Assume that $|A|\geq 3$ and let $u\in A\backslash \{0, 1\}$.
		Then we immediately have $ d_{u} $, $\chi^{(0)}$, $ \chi^{(u)} \in \Der(A)$ by Corollary \ref{pro09}. Furthermore, it is easy to see that $d_{u}\neq\textnormal{Id}_{A}$, $d_{u}\neq \textbf{0}_{A} $, $\chi^{(0)}\neq\textnormal{Id}_{A}$, $\chi^{(0)}\neq\textbf{0}_{A}$, $ \chi^{(u)}\neq\textnormal{Id}_{A}$ and $\chi^{(u)}\neq\textbf{0}_{A}$. Also, $\chi^{(0)}\neq d_{u}$, $ \chi^{(0)}\neq\chi^{(u)} $ and
		$ \chi^{(u)}\neq d_{u}$ by Lemma \ref{lem09}.
		Consequently, we have that
		$\textnormal{Id}_{A}$, $ \textbf{0}_{A} $, $ d_{u} $, $\chi^{(0)}$ and $\chi^{(u)}$ are mutually different $(\odot,\vee)$-derivations on $A$.
		\vspace{0.20cm}

		\noindent$ (2) $ Assume that $|A|\geq 4$ and let $u, v\in A\backslash \{0, 1\}$ with $u\neq v$. By Lemma \ref{lem09} (1), we have $ \chi^{(0)}\neq \textnormal{Id}_{A}$, $ \chi^{(0)}\neq \chi^{(u)}$, $ \chi^{(0)}\neq \chi^{(v)}$, $ \chi^{(0)}\neq d_{u}$ and $ \chi^{(0)}\neq d_{v}$. Clearly, $\chi^{(u)}\neq \chi^{(v)} $ and $ d_{u}\neq d_{v} $. In addition, $ d_{p}\neq \chi^{(q)} $ for any $ p,q\in \{u,v\} $ by Lemma \ref{lem09} (2). Thus we conclude that $\textnormal{Id}_{A}$, $ \textbf{0}_{A} $, $ d_{u} $, $ d_{v} $, $ \chi^{(0)} $, $ \chi^{(u)} $, $ \chi^{(v)} $ are mutually different $(\odot,\vee)$-derivations on $A$.
		\vspace{0.20cm}
		
		\noindent$ (3) $ Assume that $|A|\geq 5$. Then there exist
		$u, v\in A\backslash \{0, 1\}$ with $u<v$ (i.e, $u\leq v$ and $u\neq v$). In fact, 
		if $x, y$ are not comparable for any $x, y\in A\backslash \{0, 1\}$ with $x\neq y$, then the distributive lattice $(A, \leq)$ has a copy of $M_{5}$, which is contradicting to \cite[Theorem 3.6]{bu}.

		Let $w\in A\backslash \{0, u, v, 1\}$. By Lemma \ref{lem09} (1), we have $ \chi^{(0)}\neq \textnormal{Id}_{A}$, $ \chi^{(0)}\neq \chi^{(u)}$, $ \chi^{(0)}\neq \chi^{(v)}$, $ \chi^{(0)}\neq \chi^{(w)}$, $ \chi^{(0)}\neq d_{u}$, $ \chi^{(0)}\neq d_{v}$ and $ \chi^{(0)}\neq d_{w}$. 
		In addition, $ d_{p}\neq \chi^{(q)} $ for any $ p,q\in \{u,v,w\} $ by Lemma \ref{lem09} (2). Furthermore, $(d_{v})^{0},(d_{v})^{u}\in \Der(A) $ by Proposition \ref{pro08}. 
		By Corollary \ref{pro09}, we can get that $ (d_{v})^{r}\neq\chi^{(s)} $ for any $r\in\{0,u\}, s\in \{0,u,v,w\}$. 
		Also, $(d_{v})^{0}\neq\textbf{0}_{A}, ( d_{v})^{u}\neq d_{u} $. Indeed, if $ (d_{v})^{0}=\0 $, we have $ (d_{v})^{0}(u^{*})=v\odot u^{*}=0 $. It follows by Lemma \ref{lem01} that
		$v\leq u$, contradicting to the fact that $u<v$. If $ ( d_{v})^{u}= d_{u} $, then
		$v\odot u^{*} = (d_{v})^{u}(u^{*})= d_{u}(u^{*})= u\odot u^{*}=0$. Similarly, we get $v\leq u$, contradicting to the fact that $u<v$. 
		
		Note that $ w $ must be comparable with $ u $ or $ v $. Otherwise, if $w$ is not comparable for $u,v\in A\backslash \{0, 1\}$ with $u\leq v$, then the distributive lattice $(A, \leq)$ has a copy of $N_{5}$, which is contradicting to \cite[Theorem 3.6]{bu}. There are two cases. 
		If $ u<w $, $ v$	is not comparable for $ w $, we have $ (d_{w})^{0},(d_{w})^{u}\in \Der(A) $ and similarly, $ (d_{w})^{r}\neq\chi^{(s)} $ for any $r\in\{0,u\}, s\in \{0,u,v,w\}$. Also, it can be proved in the same way as shown before that $(d_{w})^{0}\neq0, ( d_{w})^{u}\neq d_{u} $. 
		If $ w<v $, $ u$	is not comparable for $ w $, we have $ (d_{w})^{0},(d_{v})^{w}\in \Der(A) $ and they are different from other $(\odot,\vee)$-derivations on $A$.  
		
		Finally, it is easy to check that $\textnormal{Id}_{A}$, $ \0 $, $ d_{u} $, $ d_{v} $, $ d_{w} $, $(d_{v})^{u} $, $ (d_{v})^{0} $, $(d_{w})^{0}  $, $ (d_{w})^{u}((d_{v})^{w}) $, $ \chi^{(0)} $, $ \chi^{(u)} $, $ \chi^{(v)} $, $ \chi^{(w)} $ are mutually different $(\odot,\vee)$-derivations on $A$.
	\end{proof}	
	
	\begin{proposition}\label{pro13}
		Let $A$ be an nontrivial MV-algebra. Then the following statements hold:
		\begin{enumerate}
			\item[$ (1) $] $|\Der(A)|=2$ if and only if $|A|=2$.
			
			\item[$ (2) $] $|\Der(A)|=5$ if and only if $|A|=3$.
			
			\item[$ (3) $] $|\Der(A)|=9$ if and only if $|A|=4$.	
		\end{enumerate}
	\end{proposition}
	\begin{proof}
		$(1)$ Assume that $A$ is a $2$-element MV-algebra. Then
		$A=\{0, 1\}$ is a $2$-element MV-chain, and so $|\Der(A)|=2$ by Theorem \ref{thm01}.
		
		Conversely, assume that $|\Der(A)|=2$. If $|A|\geq 3$,
		then  $|\Der(A)|\geq 5$ by Corollary \ref{cor02} $(1)$, a contradiction. Since $ A $ is nontrivial, finally we get $ |A|=2 $.
		\vspace{0.20cm}	
		
		\noindent$(2)$ Assume that $A$ is a $3$-element MV-algebra. Then  $A$ is a $3$-element MV-chain by Lemma \ref{lem12}, and so $|\Der(A)|=5$ by Theorem \ref{thm01}.
		
		Conversely, assume that $|\Der(A)|=5$. If $|A|\geq 4$, then $|\Der(A)|\geq 7$ by Corollary \ref{cor02} (2), a contradiction.
		Thus $|A|\leq 3$. But $ A $ is nontrivial and $|A|= 2$ implies that $|\Der(A)|=2$ by (1). Therefore, $|A|=3$, and consequently $A$ is a $3$-element MV-chain.
		\vspace{0.20cm}		
		
		\noindent$(3)$ Assume that $A$ is a $4$-element MV-algebra. Then  $A$ is isomorphic to the $4$-element MV-chain $L_{4}$ or the 
		$4$-element Boolean algebra $B_{4}$ by Lemma \ref{lem12}. Recall Corollary \ref{cor04} that when the MV-algebra $A$ is a Boolean algebra, $d$ is an $(\odot,\vee)$-derivation on $A$ if and only if $d$ is a derivation on the lattice $(A, \leq)$.
		It follows by Theorem \ref{thm01} and 
		\cite[Theorem 3.21]{DL} that $|\Der(A)|=9$.
		
		Conversely,  assume that $|\Der(A)|=9$. If $|A|\geq 5$, then $|\Der(A)|\geq 13$ by Corollary \ref{cor02} (3), a contradiction. Thus $|A|\leq 4$. But $ A $ is nontrivial and
		Items $(1)$ and $(2)$ imply that
		$|A|\neq 2$ and  $|A|\neq 3$. Therefore, $|A|=4$.	
	\end{proof}
	
	\subsection{Isotone $(\odot,\vee)$-derivations on MV-algebras}\label{sec3.3}
	\ 
	\newline 
	\indent
	~~~~In this subsection, we consider the condition when an $(\odot,\vee)$-derivation $ d $ is isotone and characterize the properties of the fixed point set of $ d $. 
	\begin{definition}
		\textnormal{	Let $A$ be an MV-algebra and $d\in \operatorname{Der}(A)$. $d$    is called \textbf{isotone} if for all
			$x, y \in A$, $x \leq y$ implies that $d(x) \leq d(y)$.
		}	
	\end{definition}
	
	It is clear that $ \textnormal{Id}_{A}$ and $ \0 $ are isotone. Furthermore, we have:
	
	\begin{lemma}\label{lem07}
		Let $ A $ be an MV-algebra and $a \in A$. Then the principal $(\odot,\vee)$-derivation $d_{a}$ is  isotone. 
	\end{lemma}
	\begin{proof}
		Let $ x, y \in A$ with $x \leq y$. Then $d_{a}(x)=a \odot x \leq a \odot y=d_{a}(y)$ by Lemma \ref{lem02} (4), and thus $d_{a}$ is isotone.
	\end{proof} 
	
	By \cite[Proposition 2.5]{DL}, we know that  a derivation $d$ on a bounded lattice $L$ is isotone iff $d$ is  principal. However, 
	there are  other  isotone $(\odot,\vee)$-derivations on an MV-algebra $A$ besides principal 
	$(\odot,\vee)$-derivations. 
	
	\begin{example}\label{exa05}
		Let $d=\chi^{(\frac{2}{3})}\in \operatorname{Der}(L_{4})$ (see Corollary \ref{pro09}), i.e, $d(0)=0, d(\frac{1}{3})=\frac{1}{3}, d(\frac{2}{3})=\frac{2}{3}, d(1)=\frac{2}{3} $. Then $d$
		is isotone, while $d$ is not  principal, since 
		$d(1) = \frac{2}{3}=\frac{2}{3} \odot 1$ but $  d(\frac{1}{3})=\frac{1}{3}\neq 0= \frac{2}{3}\odot \frac{1}{3}$. 
	\end{example}
	
	Proposition \ref{pro04} says that
	if $d$ is an $(\odot,\vee)$-derivation on an MV-algebra $A$ with $d(1)\in \textnormal{\textbf{B}}(A)$, then
	$d$ is isotone iff $d$ is principal.

	
	%
	%
	%
	%
	
	\begin{proposition}\label{pro04}	
		Let $A$ be an MV-algebra and $d\in \operatorname{Der}(A)$ with $d(1)\in \textnormal{\textbf{B}}(A)$. Then the following statements are equivalent: 
		\begin{enumerate}
			\item[$(1)$] $ d $ is isotone;
			
			\item[$(2)$] $d(x) \leq d(1)$ for any $x \in A$;
			
			\item[$(3)$] $d(x) = d(1) \odot x$ for any $x \in A$;
			
			\item[$(4)$] $d(x \wedge y) = d(x) \wedge d(y)$ for all $x, y \in A$;
			
			\item[$(5)$] $d(x \vee y) = d(x) \vee d(y)$ for all $x, y \in A$;
			
		\end{enumerate}
	\end{proposition}
	
	\begin{proof}
		
		\noindent$(1)\Rightarrow (2)$ is clear since $ x\leq 1 $ holds for any $ x\in A $.
		
		\noindent$(2)\Rightarrow (3).$ Assume that $d(x) \leq d(1)$ for any $x \in A$.
		Since $d(x) \leq x$ by Proposition \ref{pro01} (4), it follows that
		$$d(x) \leq d(1)\wedge x=d(1)\odot x \leq d(x) $$
		by Lemma \ref{lem03} (6) and Proposition \ref{pro01} (5). Thus
		$d(x)=d(1)\odot x$.
		\vspace{0.20cm}
		
		\noindent$(3)\Rightarrow (4).$ 
		Assume that $d(x) = d(1) \odot x$ for any $x \in A$. Then for all $x, y\in A$, we have
		$d(x\wedge y)=d(1)\odot(x\wedge y)=(d(1)\odot x)\wedge (d(1)\odot y)=d(x)\wedge d(y)$
		by Lemma \ref{lem02} (6).
		\vspace{0.20cm}
		
		\noindent$(4)\Rightarrow (1).$ Assume that $d(x \wedge y) = d(x) \wedge d(y)$ for all $x, y \in A$.
		Let
		$x\leq y$. Then $ d(x)=d(x\wedge y)=d(x) \wedge d(y) \leq d(y) $, and thus $d$ is isotone.
		\vspace{0.20cm}
		
		\noindent$(3)\Rightarrow (5).$
		Assume that $(3)$ holds. Then for all $x, y\in A$, we have
		$d(x\vee y)=d(1)\odot(x\vee y)=(d(1)\odot x)\vee (d(1)\odot y)=d(x)\vee d(y)$	
		by	Lemma \ref{lem02} (7).	
		\vspace{0.20cm}
		
		\noindent$(5)\Rightarrow (1).$ 
		Assume that $(5)$ holds. Then for all $x, y\in A$ with $x\leq y$, we have $ d(x)\leq d(x)\vee d(y)=d(x\vee y)=d(y)$, and thus $d$ is isotone.	
		
	\end{proof}

	\begin{corollary}\label{cor03}
		Let $ A $ be an MV-algebra. 
		Denote the set of all isotone $(\odot,\vee)$-derivations with $d(1)\in \textnormal{\textbf{B}}(A) $
		by $ \IDer(A)$, i.e,
		$ \IDer(A)=\{ d\in \Der(A)~ |~ d$ is isotone and $d(1)\in \textnormal{\textbf{B}}(A)\}$. Then there is a bijection between $\IDer(A)$ and $ \textnormal{\textbf{B}}(A) $.
	\end{corollary}
	\begin{proof}
		Define a map $ f : \IDer(A) \rightarrow \textbf{B}(A) $ by $f(d) = d(1) $ for any $ d \in \IDer(A) $. And define a map
		$ g : \textbf{B}(A)  \rightarrow \IDer(A) $ by $ g(a) = d_{a} $ for any $ a \in A $. Then by Proposition \ref{pro04}, we have $ fg = \textnormal{Id}_{A} $ and
		$ gf = \textnormal{Id}_{\textnormal{IDer(A)}} $. Hence $ f $ is a bijection.
	\end{proof}
	
	Generally, $d$ is an isotone $(\odot, \vee)$-derivation on an MV-algebra $A$ does not necessarily imply that $d(x \oplus y) = d(x) \oplus d(y)$ for all $x, y \in A$. For example, In the MV-algebra $L_{3}$, $ \chi^{(\frac{1}{2})}\in \Der(A) $ and is isotone while $ \chi^{(\frac{1}{2})}(\frac{1}{2}\oplus\frac{1}{2})=\chi^{(\frac{1}{2})}(1)=\frac{1}{2}\neq 1=\frac{1}{2}\oplus \frac{1}{2}=\chi^{(\frac{1}{2})}(\frac{1}{2})\oplus \chi^{(\frac{1}{2})}(\frac{1}{2}) $. In the following proposition, the condition $d(1)\in \textnormal{\textbf{B}}(A)$ cannot be removed. 
	
	\begin{proposition}\label{pro:3.20}
		Let $A$ be an MV-algebra, and $d\in \operatorname{Der}(A)$.  Then the following statements are equivalent: 
		\begin{enumerate}
			\item[$(1)$] $d\in \IDer(A)$;
			\item[$(2)$] $d(x \oplus y) = d(x) \oplus d(y)$ for all $x, y \in A$;
			\item[$(3)$] $d(x \odot y) = d(x) \odot d(y)$ for all $x, y \in A$.
		\end{enumerate}
	\end{proposition}
	
	\begin{proof}
		\noindent$(1)\Rightarrow (2).$ Assume $ d\in\IDer(A) $. By Lemma \ref{lem14}, supposing that $ A $ is a subdirect product of a family $ \{A_{i}\}_{i\in I} $ of MV-chains, let $ h: A \rightarrow \prod_{i\in I}A_{i}$ be a one-one homomorphism and  for each $ j\in I $, the composite map $ \pi_{j}\circ h $ is a 
		homomorphism onto $ A_{j}$. Let $d(1)=a=(a_{i})_{i\in I}\in \textnormal{\textbf{B}}(A)$. Then $a_{i}\in \textnormal{\textbf{B}}(A_{i})$ and by Lemma \ref{lem15} (1) we have $ a_{i}=0$ or $a_{i}=1$ for each $ i\in I $. Since $ d\in \IDer(A) $, it follows by Proposition \ref{pro04} that for any $ x=(x_{i})_{i\in I}\in A $, $ d(x)=x\odot a=(x_{i}\odot a_{i})_{i\in I} $. Therefore, $ d(x\oplus y)=((x_{i}\oplus y_{i})\odot_{i} a_{i})_{i\in I} =((x_{i}\odot_{i} a_{i})\oplus(y_{i}\odot_{i} a_{i}))_{i\in I}=d(x)\oplus d(y)$.
		
		\vspace{0.20cm}   
		\noindent$(2)\Rightarrow (1).$  Assume that $d(x \oplus y) = d(x) \oplus d(y)$ for all $x, y \in A$, we immediately get $ d(1)=d(1)\oplus d(1) $. Hence $d(1)\in \textnormal{\textbf{B}}(A)$. To prove that $d$ is isotone, let
		$x\leq y$. Then by Lemma \ref{lem01} (4) there exists an element $ z\in A $ such that $y=x\oplus z $. So $ d(y)=d(x\oplus z)=d(x) \oplus d(z) \leq d(x) $, and thus $d$ is isotone.
		\vspace{0.15cm}
		
		\noindent$(1)\Rightarrow (3).$ 
		Assume that $(1)$ holds, by Proposition \ref{pro04} we have $ d(x)=d(1)\odot x $.	Since $d(1)\in \textnormal{\textbf{B}}(A)$,
		$d(x\odot y)=d(1)\odot(x\odot y)=(d(1)\odot x)\odot (d(1)\odot y)=d(x)\odot d(y)$ for all $x, y\in A$.
		\vspace{0.15cm}
		
		\noindent$(3)\Rightarrow (1).$ Assume that $(3)$ holds. Then for any $x\in A$, we have 
		$ d(x)=d(x\odot 1)=d(x)\odot d(1)\leq d(1) $ 
		by Lemma \ref{lem02} (1). Set $ x=y=1 $ in $ (3) $, we have $ d(1)=d(1)\odot d(1) $. Hence $d(1)\in \textnormal{\textbf{B}}(A)$. 
		Thus by Propostion \ref{pro04} we get $ d$ is isotone. Therefore, $ (1) $ holds.
	\end{proof} 
	
	%
	%
	%
	%
	
	\begin{corollary}\label{cor 001}	
		Let $A$ be an MV-algebra and $d\in \IDer(A)$. Then	
		$ d $ is idempotent, that is, $d^{2}=d$. 
	\end{corollary}
	\begin{proof}
		Assume that $d\in \IDer(A)$. By Proposition \ref{pro04} (3) and Proposition \ref{pro:3.20}, we have $ d(d(x))=d(1\odot d(x))=d(1)\odot d(x)=d(1\odot x)= d(x) $ for any $x\in A$. Thus $d^{2}=d$.
	\end{proof}
	
	Generally, the converse of Corollary \ref{cor 001} does not hold. For example, let
	$ d=\chi^{(0)} \in \Der (L_{3}) $. Then $ d(1)=0\in \textnormal{\textbf{B}}(A) $ and $ d $ is idempotent. But $ d $ is not isotone, since $ d(\frac{1}{2})=\frac{1}{2}> 0=  d(1) $.

	Using the fixed point sets of isotone derivations, the characterizations of some different types of lattice have been described in \cite{lattice2012}. Analogously, we next discuss the relation between ideals and fixed point sets of $(\odot,\vee)$-derivations on MV-algebras.
	
	Let $A$ be an MV-algebra and $d\in \operatorname{Der}(A)$. Denote the set of all \textbf{fixed point of $d$} by $ \F_{d}(A) $, i.e.,
	$$ \F_{d}(A)=\{x\in A~|~d(x)=x\}. $$ 
	By Proposition \ref{pro01} (10),
	$ \F_{d}(A) $ is a downset.
	
	\begin{proposition}\label{pro10}
		Let $A$ be an MV-algebra. If $d\in \PDer(A)$, then $ \F_{d}(A) $ is a lattice ideal of $ A $. 
	\end{proposition}
	\begin{proof}	
		Assume that $d\in \PDer(A)$, and let $d=d_{a}$, where $a\in A$. Then $d(x)=a\odot x$ for any $x\in A$.
		To prove that $ \F_{d}(A) $ is closed under $\vee$, let  $x, y \in \F_{d}(A) $. Then $d(x) = x $ and $d(y) = y$. It follows  by Lemma \ref{lem02} (7) that $ d(x \vee y) =a\odot(x\vee y)=(a\odot x)\vee(a\odot y)=d(x) \vee d(y) = x \vee y$, and so $x \vee y \in \F_{d}(A)$. Thus $ \F_{d}(A) $ is closed under $\vee$. This, together with the fact that $ \F_{d}(A) $ is a downset,
		implies that $ \F_{d}(A) $ is a lattice ideal of $ A $. 
	\end{proof} 
	
	
	
	
	\section{Direct product of  $(\odot,\vee)$-derivations}\label{s04}
	~~~~In this section, we will discuss the relation between direct product of  $(\odot,\vee)$-derivations and $(\odot,\vee)$-derivations on the direct product of MV-algebras.
	
	\begin{definition}\cite{MV4}
		\textnormal{Let $\Omega$ be an index set. The \textbf{direct product} $\prod_{i \in \Omega} A_{i}$ of a family $\left\{A_{i}\right\}_{i \in \Omega}$ of MV-algebras is the MV-algebra obtained by endowing the set-theoretical cartesian product of the family with the MV-operations defined pointwise. In other words, $\prod_{i \in \Omega} A_{i}$ is the set of all functions $f: \Omega \rightarrow \bigcup_{i \in \Omega} A_{i}$ such that $f(i) \in A_{i}$ for all $i \in \Omega$, with the operations `` $*$ '' and `` $\oplus$ '' defined by
			$$
			(f^{*})(i)={ }_{\text {def }}  (f(i))^{*} \quad\text { and } \quad (f \oplus g)(i)=_{\text{def}} f(i) \oplus g(i)  ~\text{for ~all}~ i\in \Omega.
			$$
			The zero element $ 0 $ of $\prod_{i \in \Omega} A_{i}$ is the function $i \in \Omega \mapsto 0_{i} \in A_{i}$, and  the element $ 1 $ of $\prod_{i \in \Omega} A_{i}$ is the function $i \in \Omega \mapsto 1_{i} \in A_{i}$ for all $i\in \Omega$.
		}
	\end{definition}
	
	The binary operation `` $ \odot $ " and `` $ \ominus $ " on $\prod_{i \in \Omega} A_{i}$ can be induced by `` $\oplus$ " and `` $^{*}$ ".
	Let $g, h \in \prod_{i \in \Omega} A_{i} .$ By Lemma \ref{lem01} we know that $g \leqslant h$ in $\prod_{i \in \Omega} A_{i}$ if and only if $g^{*}\oplus h=1$ if and only if 
	$(g(i))^{*}\oplus h(i)=1_{i}  $ in $A_{i}$
	if and only if
	$g(i) \leqslant h(i)$  for any $i \in \Omega$. As usual, we write $(g(i))_{i \in \Omega}$ for $g$.
	
	\begin{definition}\cite{MV4}
		\textnormal{For each $i \in \Omega$, define the map $\pi_{i}: \prod_{i \in \Omega} A_{i} \rightarrow A_{i}$ by $\pi_{i}(g)=g(i)$ for any $g \in \prod_{i \in \Omega} A_{i}$, and define the map $\rho_{i}: A_{i} \rightarrow \prod_{i \in \Omega} A_{i}$ by
			$$
			\left(\rho_{i}(a)\right)(j)= \begin{cases}a, & \text { if } j=i \\ 0_{j}, & \text { otherwise }\end{cases}
			$$
			for any $a \in A_{i}$. $\pi_{i}$ is called the \textbf{$i$-th projection}, and $\rho_{i}$ is called the \textbf{$i$-th embedding}.
		}
	\end{definition}
	\begin{definition}
		\textnormal{
			For each $i \in \Omega$, let $d_{i}$ be an operator on $A_{i} .$ Define an operator $\prod_{i \in \Omega} d_{i}: \prod_{i \in \Omega} A_{i} \rightarrow \prod_{i \in \Omega} A_{i}$ by $\left(\prod_{i \in \Omega} d_{i}\right)(g)=\left(d_{i}(g(i))\right)_{i \in \Omega}$ for any $g \in \prod_{i \in \Omega} A_{i}$, and we call $\prod_{i \in \Omega} d_{i}$ the \textbf{direct product of the $\{d_{i}\}_{i \in \Omega}$}. 
		}
	\end{definition}
	When $\Omega=\{1,2, \cdots, n\}$, we denote the direct product of $\left\{A_{i}\right\}_{i \in \Omega}$ and the direct product of $\left\{d_{i}\right\}_{i \in \Omega}$, respectively, by $A_{1} \times A_{2} \times \cdots \times A_{n}$ and $d_{1} \times d_{2} \times \cdots \times d_{n}$.

	\begin{lemma}\label{lem08}
		Let $\Omega$ be an index set, $\left\{A_{i}\right\}_{i \in \Omega}$ be a family of MV-algebras, and $ d $ be an operator on $\prod_{i \in \Omega} A_{i}$. Then the following statements hold:
		\begin{enumerate}
			\item[$ (1) $] $d \in \Der(\prod_{i \in \Omega} A_{i})$ implies that $\pi_{i} d \rho_{i} \in \Der(A_{i})$ for each $i \in \Omega$;
			
			\item[$ (2) $] $d \in \Der(\prod_{i \in \Omega} A_{i})$ and $d$ is isotone implies that $\pi_{i} d \rho_{i} \in \Der(A_{i})$ and is isotone for each $i \in \Omega$;
			
			
			\item[$ (3) $] $d \in \PDer(\prod_{i \in \Omega} A_{i})$ implies that $\pi_{i} d \rho_{i} \in \PDer(A_{i})$ for each $i \in \Omega$.
		\end{enumerate}
	\end{lemma}
	\begin{proof}
		$ (1) $ Assume that $d \in \Der(\prod_{i \in \Omega} A_{i})$. For each $i \in \Omega  $, let $ x,y\in A_{i} $. Then we have 
		$$\begin{aligned} 
		\left(\pi_{i} d \rho_{i}\right)(x \odot y)&=\pi_{i} d\left(\rho_{i}(x \odot y)\right)=\pi_{i}\left(d\left(\rho_{i}(x) \odot \rho_{i}(y)\right)\right) \\
		&= \pi_{i}\left(\left(d\left(\rho_{i}(x)\right) \odot \rho_{i}(y)\right) \vee\left(\rho_{i}(x) \odot d\left(\rho_{i}(y)\right)\right)\right) \\
		&=\left(\pi_{i} d \rho_{i}(x) \odot \pi_{i} \rho_{i}(y)\right) \vee\left(\pi_{i} \rho_{i}(x) \odot \pi_{i} d \rho_{i}(y)\right) \\
		&=\left(\pi_{i} d \rho_{i}(x) \odot y\right) \vee\left(x \odot \pi_{i} d \rho_{i}(y)\right) 
		\end{aligned}$$
		and so $\pi_{i} d \rho_{i} \in \Der(A_{i})$.
		
		\vspace{0.20cm} 
		\noindent$ (2) $ Assume that $d \in \Der(\prod_{i \in \Omega} A_{i})$ and $ d $ is isotone. For each $i \in \Omega$, 
		we know  by (1) that $\pi_{i} d \rho_{i} \in \Der(A_{i})$. 
		Also, since  $\pi_{i}$ and $\rho_{i}$ are isotone, it follows that $\pi_{i} d \rho_{i}$ is isotone. Thus $(2)$ holds.
		
		
		\vspace{0.20cm} 
		\noindent$ (3) $  Assume that $d \in \PDer(\prod_{i \in \Omega} A_{i})$, i.e, $d=d_{a}$ for some $a=\left(a_{i}\right)_{i \in \Omega} \in \prod_{i \in \Omega} A_{i}$. For each $i \in \Omega$, let $x \in A_{i}$. Then we have $$\left(\pi_{i} d \rho_{i}\right)(x)=\pi_{i} d\left(\rho_{i}(x)\right)=\pi_{i}\left(\rho_{i}(x) \odot a\right)=\pi_{i}\left(\rho_{i}(x)\right) \odot \pi_{i}(a)=x \odot a_{i},$$ and thus $\pi_{i} d \rho_{i} \in \PDer(A_{i})$.
	\end{proof}
	
	Combining the structures of an MV-algebra and an $(\odot,\vee) $-derivation in the language of universal algebra~\cite{bu}, we give
	
	\begin{definition}
		A \textbf{differential MV-algebra} is an algebra $(A, \oplus, *, d, 0)$ of type $(2, 1, 1, 0)$ such that
		\begin{enumerate}
			\item[(1)] $(A, \oplus, *, 0)$ is an MV-algebra, and
			
			\item[(2)]  $d$ is an $(\odot,\vee) $-derivation on $A$.
		\end{enumerate}
	\end{definition}
	
	Let $\Omega$ be an index set, $\left\{A_{i}\right\}_{i \in \Omega}$ be a family of MV-algebras, and
	$d_{i} \in \Der(A_{i})$. Then $\left(A_{i}, \oplus_{i}, *_{i}, d_{i}, 0_{i}\right)$ is a differential MV-algebra. From the viewpoint of universal algebra \cite[Theorem 11.9]{bu}, we know that the class of all differential MV-algebras forms a variety. Thus the direct product $\left(\prod_{i \in \Omega} A_{i}, \oplus, *, \prod_{i \in \Omega} d_{i}, 0\right)$ is also a differential MV-algebra, and so $\prod_{i \in \Omega} d_{i} \in \Der(\prod_{i \in \Omega} A_{i})$. Hence we obtain that
	
	\begin{equation}\label{equ02}
	\prod_{i \in \Omega} \Der(A_{i}) \subseteq \Der(\prod_{i \in \Omega} A_{i})
	\end{equation}
	But 
	$\prod_{i \in \Omega} \Der(A_{i}) \neq \Der(\prod_{i \in \Omega} A_{i})
	$ whenever $|\Omega|\geq 2$, see Remark \ref{rem01}.

	\begin{example}\label{exa06}
		\begin{enumerate}
			\item[$ (1) $]Let $L_{2}=\{0,1\}$ be the 2-element MV-chain. Then   $ \Der(L_{2})=\{\textnormal{Id}_{L_{2}},\textbf{0}_{L_{2}}\} $ by Theorem \ref{thm01},
			%
			so $ \Der(L_{2})\times \Der(L_{2})=\{d_{1}=\textnormal{Id}_{L_{2}}\times \textnormal{Id}_{L_{2}},d_{2}=\textnormal{Id}_{L_{2}}\times \textbf{0}_{L_{2}},d_{3}=\textbf{0}_{L_{2}}\times \textnormal{Id}_{L_{2}},d_{4}=\textbf{0}_{L_{2}}\times \textbf{0}_{L_{2}} \}\subseteq \Der(L_{2}\times L_{2}) $.
			Notice that in \cite{DL}, $ L_{2} \times L_{2} $ is denoted by $ M_{4} $, and $ d_{1},d_{2},d_{3},d_{4} $ are denoted by $\textnormal{Id}_{M_{4}},y_{2},y_{4},\textbf{0}_{M_{4}}$, respectively. By \cite[Theorem 3.21]{DL}, $ |\Der(L_{2}\times L_{2})|=9 $, so $\Der(L_{2})\times \Der(L_{2})\neq \Der(L_{2} \times L_{2})$. 
			\item[$ (2) $]
			Let $L_{3}=\{0, \frac{1}{2}, 1 \} $ be the $3$-element MV-chain with $0< \frac{1}{2}< 1$. By Theorem \ref{thm01} we have $ \Der(L_{3})=\{\textnormal{Id}_{L_{3}},\textbf{0}_{L_{3}},d_{\frac{1}{2}}, \chi^{(0)}, \chi^{(\frac{1}{2})}\} $.
			Thus
			$$|\Der(L_{2})\times\Der(L_{3})|=|\Der(L_{2})|\times 
			|\Der(L_{3})|=10.$$
			
			Let $ \textbf{0}=(0,0) $, $ \textbf{a}=(0,\frac{1}{2}) $, $ \textbf{b}=(0,1) $, $ \textbf{c}=(1,0) $, $ \textbf{d}=(1,\frac{1}{2}) $ and $ \textbf{1}=(1,1) $. Then the Hasse diagram of $ L_2\times L_3 $ is given 
			below (see Figure 1). We give all elements of $ \Der(L_{2}\times L_{3}) $ in Table \ref{t2} by Python (Full details are given in Appendix I listing \ref{a1}). It can be verified that there are $ 23 $ elements (from $ d_{11}$ to $d_{33} $) in $ \Der({L_{2}\times L_{3}}) $ but not in $  \Der(L_{2})\times\Der(L_{3}) $.
			
			\begin{table}[h]
				\centering
				\caption{$\Der( L_{2}\times L_{3} )$}
				\begin{tabular}{ccccccc|ccccccc}
					\toprule
					$ x $     & \textbf{0} & \textbf{a} & \textbf{b} & \textbf{c} & \textbf{d} & \textbf{1} & $ x $     & \textbf{0} & \textbf{a} & \textbf{b} & \textbf{c} & \textbf{d} & \textbf{1} \\
					\hline
					$ d_{1}=\textnormal{Id}_{L_{2}}\times \textnormal{Id}_{L_{3}} $ & \textbf{0} & \textbf{a} & \textbf{b} & \textbf{c} & \textbf{d} & \textbf{1} & $ d_{18} $ & \textbf{0} & \textbf{a} & \textbf{b} & \textbf{c} & \textbf{d} & \textbf{b} \\
					$ d_{2}=\textnormal{Id}_{L_{2}}\times \chi^{(\frac{1}{2})} $ & \textbf{0} & \textbf{a} & \textbf{a} & \textbf{c} & \textbf{d} & \textbf{d} & $ d_{19} $ & \textbf{0} & \textbf{a} & \textbf{b} & \textbf{c} & \textbf{d} & \textbf{d} \\
					$ d_{3}=\textnormal{Id}_{L_{2}}\times d_{\frac{1}{2}} $ & \textbf{0} & \textbf{0} & \textbf{a} & \textbf{c} & \textbf{c} & \textbf{d} & $ d_{20} $ & \textbf{0} & \textbf{a} & \textbf{a} & \textbf{c} & \textbf{d} & \textbf{0} \\
					$ d_{4}=\textnormal{Id}_{L_{2}}\times \chi^{(0)} $ & \textbf{0} & \textbf{a} & \textbf{0} & \textbf{c} & \textbf{d} & \textbf{c} & $ d_{21} $ & \textbf{0} & \textbf{a} & \textbf{b} & \textbf{0} & \textbf{0} & \textbf{a} \\
					$ d_{5}=\textnormal{Id}_{L_{2}}\times \textbf{0}_{L_{3}} $ & \textbf{0} & \textbf{0} & \textbf{0} & \textbf{c} & \textbf{c} & \textbf{c} & $ d_{22} $ & \textbf{0} & \textbf{a} & \textbf{b} & \textbf{0} & \textbf{a} & \textbf{0} \\
					$ d_{6} =\textbf{0}_{L_{2}}\times\textnormal{Id}_{L_{3}} $ & \textbf{0} & \textbf{a} & \textbf{b} & \textbf{0} & \textbf{a} & \textbf{b} & $ d_{23} $ & \textbf{0} & \textbf{a} & \textbf{b} & \textbf{c} & \textbf{c} & \textbf{0} \\
					$d_{7}= \textbf{0}_{L_{2}}\times\chi^{(\frac{1}{2})} $ & \textbf{0} & \textbf{a} & \textbf{a} & \textbf{0} & \textbf{a} & \textbf{a} & $d_{24} $ & \textbf{0} & \textbf{a} & \textbf{b} & \textbf{c} & \textbf{c} & \textbf{d} \\
					$ d_{8}=\textbf{0}_{L_{2}}\times d_{\frac{1}{2}} $ & \textbf{0} & \textbf{0} & \textbf{a} & \textbf{0} & \textbf{0} & \textbf{a} & $ d_{25} $ & \textbf{0} & \textbf{a} & \textbf{b} & \textbf{c} & \textbf{d} & \textbf{a} \\
					$ d_{9}=\textbf{0}_{L_{2}}\times\chi^{(0)} $ & \textbf{0} & \textbf{a} & \textbf{0} & \textbf{0} & \textbf{a} & \textbf{0} & $ d_{26} $ & \textbf{0} & \textbf{0} & \textbf{a} & \textbf{c} & \textbf{c} & \textbf{0} \\
					$ d_{10}=\textbf{0}_{L_{2}}\times\textbf{0}_{L_{3}} $ & \textbf{0} & \textbf{0} & \textbf{0} & \textbf{0} & \textbf{0} & \textbf{0} & $ d_{27} $ & \textbf{0} & \textbf{a} & \textbf{0} & \textbf{c} & \textbf{d} & \textbf{0} \\
					$ d_{11} $ & \textbf{0} & \textbf{a} & \textbf{a} & \textbf{c} & \textbf{d} & \textbf{c} & $ d_{28} $ & \textbf{0} & \textbf{a} & \textbf{a} & \textbf{0} & \textbf{a} & \textbf{0} \\
					$ d_{12} $ & \textbf{0} & \textbf{a} & \textbf{b} & \textbf{0} & \textbf{a} & \textbf{a} & $ d_{29} $ & \textbf{0} & \textbf{a} & \textbf{b} & \textbf{0} & \textbf{0} & \textbf{0} \\
					$ d_{13} $ & \textbf{0} & \textbf{a} & \textbf{b} & \textbf{c} & \textbf{c} & \textbf{a} & $ d_{30} $ & \textbf{0} & \textbf{0} & \textbf{0} & \textbf{c} & \textbf{c} & \textbf{0} \\
					$ d_{14} $ & \textbf{0} & \textbf{a} & \textbf{b} & \textbf{c} & \textbf{c} & \textbf{c} & $ d_{31} $ & \textbf{0} & \textbf{0} & \textbf{a} & \textbf{0} & \textbf{0} & \textbf{0} \\
					$ d_{15} $ & \textbf{0} & \textbf{a} & \textbf{b} & \textbf{c} & \textbf{d} & \textbf{0} & $ d_{32} $ & \textbf{0} & \textbf{a} & \textbf{b} & \textbf{c} & \textbf{d} & \textbf{c} \\
					$ d_{16} $ & \textbf{0} & \textbf{0} & \textbf{a} & \textbf{c} & \textbf{c} & \textbf{a} & $ d_{33} $ & \textbf{0} & \textbf{a} & \textbf{a} & \textbf{c} & \textbf{d} & \textbf{a} \\
					$ d_{17} $ & \textbf{0} & \textbf{0} & \textbf{a} & \textbf{c} & \textbf{c} & \textbf{c} & --    & --    & --    & --    & --    & --    & -- \\
					\bottomrule
				\end{tabular}%
				\label{t2}%
			\end{table}%
			\vspace{1.2cm}
			\begin{center}
				\setlength{\unitlength}{0.5cm}
				\begin{picture}(4,5)
				\thicklines
				\put(2.0,0.0){\line(-1,1){2.0}}
				\put(2.0,0.0){\line(1,1){4.0}}
				\put(4.0,6.0){\line(1,-1){2.0}}
				\put(4.0,6.0){\line(-1,-1){4.0}}
				\put(2.0,4.0){\line(1,-1){2.0}}
				\put(2.0,-1){$\textbf{0}$}
				\put(4.3,1.5){$\textbf{a}$}
				\put(6.3,3.5){$\textbf{b}$}
				\put(-1.3,2.0){$\textbf{c}$}
				\put(1.3,4.5){$\textbf{d}$}
				\put(4.0,6.5){$\textbf{1}$}
				\put(1.8,-0.2){$\bullet$}
				\put(3.8,1.8){$\bullet$}
				\put(5.8,3.8){$\bullet$}
				\put(-0.2,1.8){$\bullet$}
				\put(1.8,3.8){$\bullet$}
				\put(3.8,5.8){$\bullet$}
				\put(0.0,-2.5){Figure 1: $L_{2}\times L_{3}  $}
				\end{picture}
				
			\end{center}		
		\end{enumerate}
	\end{example}
	\vspace{1.0cm}
	
	\begin{theorem}\label{thm02}
		Let $\Omega$ be an index set, $\left\{A_{i}\right\}_{i \in \Omega}$ be a family of MV-algebras, and $ d_{i} $ be an operator on $A_{i}$ for each $ i \in \Omega $. Let $ A=\prod_{i \in \Omega} A_{i} $. Then the following statements hold:
		\begin{enumerate}
			\item[$ (1) $] $\pi_{i}\left(\prod_{i \in \Omega} d_{i}\right) \rho_{i}=d_{i}$, and $\pi_{i}\left(\prod_{i \in \Omega} d_{i}\right)=d_{i} \pi_{i}$ for each $i \in \Omega$.
			
			\item[$ (2) $] $\prod_{i \in \Omega} d_{i} \in \Der(A)$ if and only if $d_{i} \in \Der(A_{i})$ for each $i \in \Omega$.

			\item[$ (3) $] $\prod_{i \in \Omega} d_{i} \in \Der(A)$ and $\prod_{i \in \Omega} d_{i}$ is isotone  if and only if $d_{i} \in \Der(A_{i})$ and $d_{i}$ is isotone for each $i \in \Omega$.
			
			\item[$ (4) $] $\prod_{i \in \Omega} d_{i} \in \PDer(A)$ if and only if $d_{i} \in \PDer(A_{i})$ for each $i \in \Omega$.
			

			\item[$ (5) $] For any $i \in \Omega$, if $d_{i}\left(0_{i}\right)=0_{i}$, then $\left(\prod_{i \in \Omega} d_{i}\right) \rho_{i}=\rho_{i} d_{i}$, that is, the corresponding diagram is commutative (put $d=\prod_{i \in \Omega} d_{i}$).
			\begin{displaymath}
			\xymatrix{
				A_{i} \ar[r]^{d_{i}} \ar[d] _{\rho_{i}} & A_{i} \ar[d]^{\rho_{i}} \\
				\prod_{i \in \Omega} A_{i} \ar[r]_{d} & \prod_{i \in \Omega} A_{i} 
			}
			\end{displaymath}
			%
		\end{enumerate}
	\end{theorem}
	\begin{proof}
		$ (1) $ Let $i \in \Omega$ and $a \in A_{i}$. It is easy to see that $\left(\pi_{i}\left(\prod_{i \in \Omega} d_{i}\right) \rho_{i}\right)(a)=d_{i}(a)$, and so $\pi_{i}\left(\prod_{i \in \Omega} d_{i}\right) \rho_{i}=d_{i} .$ Also, for any $z \in A$, we have $\left(\pi_{i}\left(\prod_{i \in \Omega} d_{i}\right)\right)(z)=d_{i} \pi_{i}(z)$, since $z=\left(\pi_{i}(z)\right)_{i \in \Omega} .$ Thus $\pi_{i}\left(\prod_{i \in \Omega} d_{i}\right)=d_{i} \pi_{i}$.
		
		\vspace{0.20cm} 
		\noindent$ (2) $ Assume that $d_{i} \in \Der(A_{i})$ for each $i \in \Omega$. Then   $\prod_{i \in \Omega} d_{i} \in \Der(A)$ by Eq.~\eqref{equ02}. 
		
		Conversely, if $\prod_{i \in \Omega} d_{i} \in \Der(A)$, then $d_{i}=\pi_{i}\left(\prod_{i \in \Omega} d_{i}\right) \rho_{i} \in \Der(A_{i})$ by (1) and Lemma \ref{lem08} (1).
		
		\vspace{0.20cm} 
		\noindent$ (3) $ Assume that $d_{i} \in \Der(A_{i})$ for each $i \in \Omega$ and $ d_{i} $ is isotone for each $ i\in\Omega $. Then $\prod_{i \in \Omega} d_{i} \in \Der(A)$ by (2). And it can be verified that $\prod_{i \in \Omega} d_{i}$ is isotone. In fact, let $ x,y\in A $ and $ x\leq y $, that is, $x_{i}\leq y_{i}$ for each $ i\in\Omega $, we have $(\prod_{i \in \Omega} d_{i})(x)=\prod_{i \in \Omega} d_{i}(x_{i})\leq \prod_{i \in \Omega} d_{i}(y_{i})=(\prod_{i \in \Omega} d_{i})(y)$.
		
		Conversely, if $\prod_{i \in \Omega} d_{i} \in \Der(A)$ and $\prod_{i \in \Omega} d_{i}$ is isotone, then $d_{i}=\pi_{i}\left(\prod_{i \in \Omega} d_{i}\right) \rho_{i} \in \Der(A_{i})$ and $ d_{i} $ is isotone by (1) and Lemma \ref{lem08} (2).
		
		\vspace{0.20cm} 
		\noindent$ (4) $ Assume that  $d_{i}\in \PDer(A_{i})$ for each $i \in \Omega$. Then $d_{i}\left(x_{i}\right)=x_{i} \odot a_{i}$, where $a_{i} \in A_{i}$. Let $\left(x_{i}\right)_{i \in \Omega} \in A$, and so $\left(\prod_{i \in \Omega} d_{i}\right)\left(\left(x_{i}\right)_{i \in \Omega}\right)=\left(d_{i}\left(x_{i}\right)\right)_{i \in \Omega}=\left(x_{i} \odot a_{i}\right)_{i \in \Omega}=\left(x_{i}\right)_{i \in \Omega} \odot\left(a_{i}\right)_{i \in \Omega} .$ Thus $\prod_{i \in \Omega} d_{i} \in \PDer(A)$.
		
		Conversely, if $\prod_{i \in \Omega} d_{i} \in \PDer(A)$, then $d_{i}=\pi_{i}\left(\prod_{i \in \Omega} d_{i}\right) \rho_{i} \in \PDer(A_{i})$ by (1) and Lemma \ref{lem08} (3).
		
		\vspace{0.20cm} 
		\noindent$ (5) $ Assume that $d_{i}\left(0_{i}\right)=0_{i}$ for any $i \in \Omega$, and $\prod_{i \in \Omega} d_{i}=d$. To prove that $d \rho_{i}=\rho_{i} d_{i}$, let $x \in A_{i}$. We have $d \rho_{i}(x)=\rho_{i} d_{i}(x)$, since
		$$
		\pi_{j}\left(d \rho_{i}(x)\right)= \begin{cases}d_{i}(x), & \text { if } j=i \\ d_{j}\left(0_{j}\right), & \text { otherwise }\end{cases}
		$$
		and$$
		\pi_{j}\left(\rho_{i} d_{i}(x)\right)= \begin{cases}d_{i}(x), & \text { if } j=i \\ 0_{j}, & \text { otherwise. }\end{cases}
		$$
		Thus $d \rho_{i}=\rho_{i} d_{i}$.
		
	\end{proof}
	
	\begin{corollary}\label{cor01}
		Let $\Omega$ be an index set, $\left\{A_{i}\right\}_{i \in \Omega}$ be a family of MV-algebras, and d be an operator on $\prod_{i \in \Omega} A_{i}$. Put $ A=\prod_{i \in \Omega} A_{i} $. Then the following statements hold:
		\begin{enumerate}
			\item[$ (1) $] If $d \in \Der(A)$, then $d \in \prod_{i \in \Omega} \Der(A_{i})$ if and only if $d=\prod_{i \in \Omega} \pi_{i} d \rho_{i}$.
			
			
			\item[$ (2) $] If $d \in \PDer(A)$, then $d \in \prod_{i \in \Omega} \PDer(A_{i})$ if and only if $d=\prod_{i \in \Omega} \pi_{i} d \rho_{i}$.
		\end{enumerate}
	\end{corollary}
	\begin{proof}
		$ (1) $ Assume that $d \in \Der(A)$. Then $\pi_{i} d \rho_{i} \in \Der(A_{i})$ for each $i \in \Omega$ by Lemma \ref{lem08} (1), which implies that $d \in \prod_{i \in \Omega} \Der(A_{i})$ if $d=\prod_{i \in \Omega} \pi_{i} d \rho_{i}$.
		
		Conversely, if $d \in \prod_{i \in \Omega} \Der(A_{i})$, then $d=\prod_{i \in \Omega} d_{i}$ for some $d_{i} \in \Der(A_{i}).$ It follows by Theorem \ref{thm02} (1) that $\pi_{i} d \rho_{i}=d_{i}$, and so $d=\prod_{i \in \Omega} \pi_{i} d \rho_{i}$.
		
		%
		\vspace{0.20cm} 
		\noindent$ (2) $ Assume that $d \in \PDer(A)$. Then $\pi_{i} d \rho_{i} \in \PDer(A_{i})$ for each $i \in \Omega$ by Lemma \ref{lem08} (3), which implies that $d \in \prod_{i \in \Omega} \PDer(A_{i})$ if $d=\prod_{i \in \Omega} \pi_{i} d \rho_{i}$.
		
		Conversely, if $d \in \prod_{i \in \Omega} \PDer(A_{i})$, then $d=\prod_{i \in \Omega} d_{i}$ for some $d_{i} \in \PDer(A_{i}).$ It follows by Theorem \ref{thm02} (1) that $\pi_{i} d \rho_{i}=d_{i}$, and so $d=\prod_{i \in \Omega} \pi_{i} d \rho_{i}$.
	\end{proof}
	
	\begin{remark}\label{rem01}
		Let $\Omega$ be an index set with $|\Omega|\geq 2$, $\left\{A_{i}\right\}_{i \in \Omega}$ be a family of MV-algebras. Then		
		$\prod_{i \in \Omega} \Der(A_{i}) \neq \Der(\prod_{i \in \Omega} A_{i})
		$, since for any $a\in \prod_{i \in \Omega} A_{i}\backslash \{1\}$, we have
		$ \chi^{(a)}\in \Der(\prod_{i \in \Omega} A_{i}) $ by Corollary \ref{pro09}, but 
		$ \chi^{(a)}\not\in \prod_{i \in \Omega} \Der(A_{i})$. In fact,  for each $ i\in\Omega $, we have $\pi_{i} \chi^{(a)} \rho_{i}\in \Der(A_{i})$ by Lemma \ref{lem08} and
		$$ \pi_{i} \chi^{(a)} \rho_{i}(1_{i})=\pi_{i}(\chi^{(a)}(\rho_{i}(1_{i})))=\pi_{i}(\rho_{i}(1_{i}))=1_{i}.$$ It follows that $ \pi_{i} \chi^{(a)} \rho_{i}=\textnormal{Id}_{A_{i}}$ by Proposition \ref{pro03}, so	 $ \chi^{(a)}\neq \textnormal{Id}_{\prod_{i \in \Omega} A_{i}}=\prod_{i \in \Omega} \pi_{i} \chi^{(a)} \rho_{i}$. Thus $ \chi^{(a)}\not\in \prod_{i \in \Omega} \Der(A_{i})$
		by Corollary \ref{cor01} $(1)$, and hence $\prod_{i \in \Omega} \Der(A_{i}) \neq \Der(\prod_{i \in \Omega} A_{i})
		$.
	\end{remark}
	
	\begin{proposition}
		Let $\Omega$ be an index set, $\left\{A_{i}\right\}_{i \in \Omega}$ be a family of MV-algebras. Then $ \PDer(\prod_{i\in\Omega} A_{i})={\prod_{i\in \Omega}} \PDer(A_{i}) $.
	\end{proposition}
	\begin{proof}
		Firstly, we have $ {\prod_{i\in \Omega}} \PDer(A_{i})\subseteq \PDer(\prod_{i\in\Omega} A_{i}) $ by Theorem \ref{thm02} (4).
		
		To prove that $ \PDer(\prod_{i\in\Omega} A_{i})\subseteq {\prod_{i\in \Omega}} \PDer(A_{i}) $, let $ d\in\PDer(\prod_{i\in\Omega} A_{i})$.
		Then for any $ x=(x_{i})_{i\in \Omega}\in \prod_{i\in\Omega} A_{i} $, by Proposition \ref{pro04}, $ d(x)=x\odot a $ for some $ a=(a_{i})_{i\in \Omega}\in \prod_{i\in\Omega} A_{i} $, so $ (\prod_{i \in \Omega} \pi_{i} d \rho_{i})(x)= (\pi_{i} d \rho_{i}(x_{i}))_{i \in \Omega}=(\pi_{i}  (\rho_{i}(x_{i})\odot a))_{i \in \Omega}=(\pi_{i}  \rho_{i}(x_{i})\odot \pi_{i}(a))_{i \in \Omega} =(x_{i}\odot \pi_{i}(a))_{i \in \Omega}=(x_{i})_{i \in \Omega}\odot (a_{i})_{i \in \Omega}=x\odot a=d(x)$. It follows that $ d=\prod_{i \in \Omega} \pi_{i} d \rho_{i} $, and so $ d\in \PDer(\prod_{i\in\Omega} A_{i}) $  by Corollary \ref{cor01} (2).
	\end{proof}

	\section{Lattice structure of $(\odot,\vee)$-derivations on MV-algebras}\label{s05}
	
	~~~~Let $(A, \oplus, *, 0)$ be an MV-algebra and let $\mathrm{O}(A)$ be the set of all operators on $A$.  Define a relation $\preceq$ on $\mathrm{O}(A)$ by:
	\begin{center}
		$(\forall~ d, d^{\prime} \in \mathrm{O}(A))$ ~~ $d \preceq d^{\prime}$ if $d(x) \leq d^{\prime}(x)$ for any $x \in A$.	
	\end{center}
	
	It is easy to verify that $\preceq$ is a partial order on $\mathrm{O}(A)$ and 
	$\mathbf{0}_{A} \preceq d \preceq \mathbf{1}_{A}$ for any $d \in \mathrm{O}(A)$, where 
	$\mathbf{1}_{A}$ is defined by $\mathbf{1}_{A}(x):=1$ for any $x \in A$.
	For any $d \in \mathrm{Der}(A)$, we have $\mathbf{0}_{A} \preceq d \preceq \textnormal{Id}_{A}$ since $0 \leq d(x) \leq x$ for any $x \in A$.
	
	We also define the following binary operations on $\mathrm{O}(A)$. For $d, d^{\prime} \in \mathrm{O}(A)$, set
	\begin{equation}\label{equ04}
	\left(d \vee d^{\prime}\right)(x):=d(x) \vee d^{\prime}(x),~~ 
	\left(d \wedge d^{\prime}\right)(x):=d(x) \wedge d^{\prime}(x) 
	\end{equation}
	for any $ x \in A $.
	
	\begin{lemma} \label{lem11}
		Let $A$ be an MV-algebra. Then $(\mathrm{O}(A), \preceq, \0,\mathbf{1}_{A} )$ is a bounded lattice for which $d\vee d'$ and $d\wedge d'$ are, respectively, the least upper bound and the greatest lower bound of $d$ and $d'$.
	\end{lemma}
	\begin{proof}
		Recall that every MV-algebra induces a natural bounded lattice structure. Since the class of all lattices is a variety and $\mathrm{O}(A)$ is the direct product of $|A|$ copies of $A$, the lemma follows immediately from the usual notions of universal algebra \cite[Definition 7.8]{bu}.
	\end{proof}
	We next explore the partial order structure of the set of $(\odot,\vee)$-derivations on MV-algebras.
	
	\begin{lemma}\label{lem13}
		Let $A$ be an MV-algebra. Then $ d\vee d'\in\Der(A) $ for all $d, d'\in \Der(A)$.
	\end{lemma}
	\begin{proof}
		Let $d, d'\in \Der(A)$ and $x, y\in A$. Then we have
		\begin{eqnarray*}
			(d\vee d')(x\odot y)&=& d(x\odot y)\vee d'(x\odot y)\\
			&=& ((d(x)\odot y)\vee (x\odot d(y)))\vee ((d'(x)\odot y)\vee (x\odot d'(y))) \\
			&=& ((d(x)\odot y)\vee (d'(x)\odot y))\vee ((x\odot d(y))\vee (x\odot d'(y))) \\
			&=&((d(x)\vee d'(x))\odot y)\vee (x\odot(d(y)\vee d'(y)))\\
			&=& ((d\vee d')(x)\odot y)\vee (x\odot(d\vee d')(y))
		\end{eqnarray*}
		by Lemma \ref{lem02} (7), and so $ d\vee d'\in \Der(A) $.
	\end{proof}
	
	For $d, d'\in \Der(A)$, 
	note that the operator $d\wedge d'$ are not necessarily in $\Der(A)$ even if $A$ is a Boolean algebra, see
	\cite[Example 3.7 and Remark 4.2]{DL}.  
	
	
	\begin{proposition}\label{pro12}
		Let $A$ be an MV-algebra. 
		\begin{enumerate}
			\item [$ (1) $] If $d\wedge d'\in\Der(A)$ for all $d,d' \in \Der(A)$, then $(\Der(A), \vee, \wedge,\textbf{0}_A,\textnormal{Id}_{A})$ is a lattice. 
			
			\item[$ (2) $] If $A$ is a finite MV-algrbra, then $(\Der(A),\preceq, \textbf{0}_A,\textnormal{Id}_{A})$ is a lattice.
		\end{enumerate}
	\end{proposition}
	\begin{proof}
		$ (1) $ For $d, d'\in \Der(A)$ and $x, y\in A$, we have known $d\vee d'\in \Der(A)$. Assume that $d\wedge d'\in \Der(A)$ for all $d, d'\in \Der(A)$. 
		Then $(\Der(A), \preceq)$ is a sublattice of the lattice $(\mathrm{O}(A), \preceq)$ by Lemma \ref{lem11}. Thus we complete the proof.
		
		\vspace{0.20cm} 
		\noindent$ (2) $ Assume that $A$ is a finite MV-algebra, by Lemma \ref{lem13} we have $d\vee d'\in \Der(A)$ for all $d,d'\in \Der(A)$.
		Since $\Der(A)$ is finite as a subset of the finite set $\mathrm{O}(A)$, it follows that $\bigvee B\colon=\bigvee_{b\in B}b$ exists for every subset $B$ of $\Der(A)$. Noticing that $ \bigvee\emptyset=\textbf{0}_A $, hence $(\Der(A),\preceq, \textbf{0}_L,\textnormal{Id}_{A})$ is a lattice by \cite[Theorem 4.2]{bu}. 
	\end{proof}
	
	In what follows, we will describe the lattice $\Der( L_{n})~ ( n\geq2) $.
	
	\begin{lemma}\label{lem5.4}
		Let $(L, \leq)  $ be a chain with the bottom element $0$, and let $$\mathcal{A}(L)=\{(x,y)\in L\times L ~| ~y\leq x\}\backslash\{(0,0)\}. $$ Then
		$(\mathcal{A}(L), \prec )$ is a sublattice of the lattice $(L\times L, \prec)$, where $\prec $ is defined by:
		for any $ (x_{1}, y_{1}), (x_{2},y_{2}) \in L\times L$,
		\begin{center}
			$ (x_{1}, y_{1})\prec(x_{2},y_{2}) $ if and only if $ x_{1}\leq x_{2} $ and $ y_{1}\leq y_{2} $.	
		\end{center} 
	\end{lemma}
	
	\begin{proof}
		It is well known that $(L\times L, \prec)$ is a lattice and 
		for any $ (x_{1}, y_{1}), (x_{2},y_{2}) \in L\times L$,
		$$(x_{1} \vee x_{2}, y_{1} \vee y_{2})= (x_{1}, y_{1})\vee (x_{2}, y_{2}), \quad (x_{1} \wedge x_{2}, y_{1} \wedge y_{2})= (x_{1}, y_{1})\wedge (x_{2}, y_{2}).$$
		
		To prove that $(\mathcal{A}(L), \prec )$ is a sublattice of the lattice $(L\times L, \prec)$, let 
		$(a, b), (c, d)\in \mathcal{A}(L)$. Then
		$b\leq a$, $d\leq c$ and $(a, b)\neq (0, 0)$, $(c, d)\neq (0, 0)$. It follows  that $b\vee d \leq a\vee c$,
		$b\wedge d \leq a\wedge c$,
		$a\neq 0$ and $c\neq 0$, so $a\vee c\neq 0$ and $a\wedge c\neq 0$
		since $L$ is a chain.
		Thus
		$(a\vee c, b\vee d)\neq (0, 0)$, and
		$(a\wedge c, b\wedge d)\neq (0, 0)$. So 
		$(a, b)\vee (c, d)\in \mathcal{A}(L) $
		and $(a, b)\wedge (c, d)\in \mathcal{A}(L) $. Consequently, we get that $(\mathcal{A}(L), \prec )$ is a sublattice of the lattice $(L\times L, \prec)$.
	\end{proof}
	
	\begin{lemma}\label{lem5.5}
		Let  $ n\geq2 $ be a positive integer,  $L_{n}$ be the $n$-element MV-chain,  and let $\mathcal{A}(L_{n})=\{(x,y)\in L_{n}\times L_{n} ~|y\leq x\}\backslash\{(0,0)\} $. Then
		the following statements hold:	
		
		\begin{enumerate}
			\item [$ (1) $] $(d_{x})^{y}\neq (d_{z})^{w}$ for any $(x, y), (z, w)\in  \mathcal{A}(L_{n})$ with
			$(x, y)\neq (z, w)$,  where $(d_{x})^{y}$ is defined by
			$$(d_{x})^{y}(z):=
			\begin{cases}
			y& \textrm{if} ~~ z = 1\\
			d_{x}(z)=x\odot z& \textrm{otherwise}.
			\end{cases}$$	 
			(See Proposition \ref{pro08}).
			
			\item[$ (2) $] $ \Der(L_{n})= \{(d_{x})^{y}~|~(x, y)\in  \mathcal{A}(L_{n}) \}.$	
			
			\item[$ (3) $]		
			$(d_{x})^{y}\wedge( d_{z})^{w}=(d_{x\wedge z})^{y\wedge w}$ and $(d_{x})^{y}\vee( d_{z})^{w}=(d_{x\vee z})^{y\vee w}$		
			for any $ (x,y),(z,w)\in\mathcal{A}(L_{n}) $.

			\item[$ (4) $] $\Der(L_{n})$ is a sublattice of $(\operatorname{O}(L_{n}), \prec)$.
		\end{enumerate}

	\end{lemma}
	\begin{proof}
		
		\vspace{0.10cm} 
		\noindent$ (1) $ Let $(x, y), (z, w)\in  \mathcal{A}(L_{n})$ with
		$(x, y)\neq (z, w)$. Then $y\leq x$, $x\neq 0$, and $w\leq z$, $z\neq 0$. So
		$x^{*}\neq 1$ and $z^{*}\neq 1$.
		
		If $y\neq w$, then   $(d_{x})^{y}(1)=y\neq w= (d_{z})^{w}(1)$, and so $(d_{x})^{y}\neq (d_{z})^{w}$. 
		
		If
		$x\neq z$ and $y=w$, then we also have 
		$(d_{x})^{y}\neq (d_{z})^{w}$. Indeed, suppose on the contrary that $(d_{x})^{y}= (d_{z})^{w}$. 
		Since $x^{*}\neq 1$ and $z^{*}\neq 1$, we have
		$$z\odot x^{*}=(d_{z})^{w}(x^{*})=(d_{x})^{y}(x^{*})=x \odot x^{*}=0 $$ and 
		$x \odot z^{*}=(d_{x})^{y}(z^{*})=(d_{z})^{w}(z^{*})=z\odot z^{*}=0 $, which implies that $z\leq x$ and $x\leq z$ by Lemma \ref{lem01}, and so $ x=z $, a contradiction.
		
		\vspace{0.20cm} 
		\noindent$ (2) $ Denote the set $\{(d_{x})^{y}~|~(x, y)\in  \mathcal{A}(L_{n}) \}$ by $\mathcal{B}$. For any $(x, y)\in  \mathcal{A}(L_{n})$, we have $y\leq x= d_{x}(1)$ and so 	$ (d_{x})^{y}\in \Der(L_{n})$ by Proposition \ref{pro08}. Thus	$ \mathcal{B}\subseteq \Der(L_{n})$.
		Also, by Item (1) we obtain that $ |\mathcal{B}|=|\mathcal{A}(L_{n})|=\frac{n(n+1)}{2}-1=\frac{(n+2)(n-1)}{2} $, so $|\mathcal{B}|=|\Der(L_{n})|$
		by Theorem \ref{thm01}. Hence $\mathcal{B}=\Der(L_{n})$.	
		
		\vspace{0.20cm} 
		\noindent$ (3) $ Let $ (x,y),(z,w)\in\mathcal{A}(L_{n})$. Then $ ((d_{x})^{y}\wedge( d_{z})^{w})(1)=(d_{x})^{y}(1)\wedge( d_{z})^{w}(1)=y\wedge w=(d_{x\wedge z})^{y\wedge w}(1)$ and  $((d_{x})^{y}\vee( d_{z})^{w})(1)=(d_{x})^{y}(1)\vee( d_{z})^{w}(1)=y\vee w=(d_{x\vee z})^{y\vee w}(1)$.
		
		Also, for $ c\in L_{n}\backslash \{1\} $, we have $$ ((d_{x})^{y}\wedge( d_{z})^{w})(c)=(d_{x})^{y}(c)\wedge( d_{z})^{w}(c)=(x\odot c)\wedge(z\odot c)=(x\wedge z)\odot c=(d_{x\wedge z})^{y\wedge w}(c) $$ by Lemma \ref{lem02} (6), and
		$$((d_{x})^{y}\vee( d_{z})^{w})(c)=(d_{x})^{y}(c)\vee( d_{z})^{w}(c)=(x\odot c)\vee(z\odot c)=(x\vee z)\odot c=(d_{x\vee z})^{y\vee w}(c)$$ by Lemma \ref{lem02} (7). 
		It follows that $ (d_{x})^{y}\wedge( d_{z})^{w}=(d_{x\wedge z})^{y\wedge w} $ and  $(d_{x})^{y}\vee( d_{z})^{w}=(d_{x\vee z})^{y\vee w}$. 
		
		\vspace{0.20cm} 
		\noindent$ (4) $	It follows immediately by Items $(2)$, $(3)$ and Lemma \ref{lem5.4} that $\Der (L_{n})$ is closed under $\vee$ and $\wedge$, so  $\Der(L_{n})$ is a sublattice of $(\operatorname{O}(L_{n}), \preceq)$.
	\end{proof}
	
	\begin{theorem}\label{thm04}
		Let  $ n\geq2 $ be a positive integer,
		$L_{n}$ be the $n$-element MV-chain,  and let
		$\mathcal{A}(L_{n})=\{(x,y)\in L_{n}\times L_{n} ~|y\leq x\}\backslash\{(0,0)\} $. Then the lattice $ \Der(L_{n}) $ is isomorphic to the lattice $\mathcal{A}(L_{n})$(see the following diagram). 
	\end{theorem}
	\begin{figure}[htbp]
		\centering
		\includegraphics[width=0.6\textwidth]{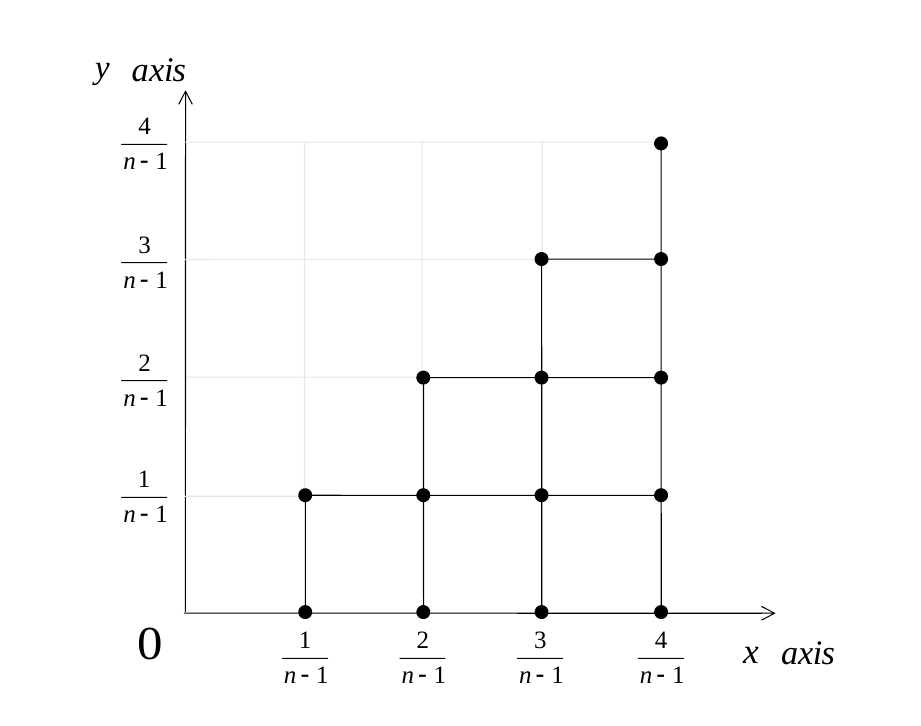}
	\end{figure}
	\begin{proof}
		Let $\mathcal{B}=\{(d_{x})^{y}|(x,y)\in\mathcal{A}(L_{n})\} $. Then $ \mathcal{B}=\Der(L_{n})$ by Lemma \ref{lem5.5} (2).		 
		
		Define a map $f\colon \mathcal{A}(L_{n})\rightarrow\mathcal{B}$
		by $(x, y)\mapsto (d_{x})^{y}$ for any $(x, y)\in \mathcal{A}(L_{n})$. Then $f$ is injective by Lemma \ref{lem5.5} $(1)$. Also, it is clear that $f$ is surjective by the definition of $ \mathcal{B} $. 
		
		To prove that $ f $ is a homomorphism, let $ (x,y),(z,w)\in\mathcal{A}(L_{n}) $. Then, by Lemma \ref{lem5.5}, we have $ f((x, y)\vee (z, w))=
		f((x\vee z, y\vee w))=(d_{x\vee z})^{y\vee w}=(d_{x})^{y}\vee( d_{z})^{w}=f((x, y))\vee f((z, w))$ and
		$ f((x, y)\wedge (z, w))=
		f((x\wedge z, y\wedge w))=(d_{x\wedge z})^{y\wedge w}=(d_{x})^{y}\wedge( d_{z})^{w}=f((x, y))\wedge f((z, w))$. Thus $f$ is a lattice isomorphism.
	\end{proof}

	\begin{example}\label{exa07}
		\begin{enumerate}
			\item [$ (1) $] We draw Hasse diagrams of $ \Der(L_{n})(2\leq n\leq 5) $ in the following:
			
			\setlength{\unitlength}{1cm}
			\begin{picture}(2,6)
			\thicklines
			\put(1.0,2.0){\line(0,1){1.0}}
			\put(0.9,2.94){$\bullet$}
			\put(0.9,1.94){$\bullet$}
			\put(0.3,0.0){$ \Der(L_{2}) $}
			\end{picture}
			\setlength{\unitlength}{0.8cm}
			\begin{picture}(4,6)
			\thicklines
			\put(1.0,3.0){\line(1,-1){1.0}}
			\put(1.0,3.0){\line(1,1){1.0}}
			\put(3.0,3.0){\line(-1,-1){1.0}}
			\put(3.0,3.0){\line(-1,1){1.0}}
			\put(2.0,5.0){\line(0,-1){1.0}}
			\put(1.88,1.9){$\bullet$}
			\put(0.88,2.9){$\bullet$}
			\put(2.88,2.9){$\bullet$}
			\put(1.88,3.9){$\bullet$}
			\put(1.88,4.96){$\bullet$}
			\put(1.3,0.0){$\Der(L_{3})$}
			\end{picture}
			\setlength{\unitlength}{0.7cm}
			\begin{picture}(5,5)
			\thicklines
			\put(2.0,2.0){\line(1,1){2.0}}
			\put(2.0,2.0){\line(-1,1){1.0}}
			\put(1.0,3.0){\line(1,1){2.0}}
			\put(1.0,5.0){\line(1,-1){2.0}}
			\put(1.0,5.0){\line(1,1){1.0}}
			\put(2.0,6.0){\line(1,-1){2.0}}
			\put(2.0,6.0){\line(0,1){1.3}}
			\put(1.8,1.9){$\bullet$}
			\put(2.8,2.85){$\bullet$}
			\put(3.86,3.85){$\bullet$}
			\put(0.8,2.9){$\bullet$}
			\put(1.84,3.87){$\bullet$}
			\put(2.84,4.87){$\bullet$}
			\put(0.86,4.9){$\bullet$}
			\put(1.84,5.9){$\bullet$}
			\put(1.86,7.2){$\bullet$}
			\put(1.4,0.0){$ \Der(L_{4}) $}
			\end{picture}
			\setlength{\unitlength}{0.6cm}
			\begin{picture}(5,5)
			\thicklines
			\put(2.0,2.0){\line(1,1){3.0}}
			\put(2.0,2.0){\line(-1,1){1.0}}
			\put(1.0,3.0){\line(1,1){3.0}}
			\put(1.0,7.0){\line(1,-1){3.0}}
			\put(1.0,7.0){\line(1,1){1.0}}
			\put(2.0,8.0){\line(1,-1){3.0}}
			\put(1.0,5.0){\line(1,1){2.0}}
			\put(1.0,5.0){\line(1,-1){2.0}}
			\put(2.0,8.0){\line(0,1){1.3}}
			\put(1.8,1.85){$\bullet$}
			\put(2.8,2.85){$\bullet$}
			\put(3.8,3.85){$\bullet$}
			\put(4.8,4.85){$\bullet$}
			\put(0.8,2.9){$\bullet$}
			\put(1.8,3.85){$\bullet$}
			\put(2.8,4.85){$\bullet$}
			\put(3.8,5.85){$\bullet$}
			\put(0.8,4.85){$\bullet$}
			\put(1.8,5.85){$\bullet$}
			\put(2.8,6.85){$\bullet$}
			\put(0.8,6.85){$\bullet$}
			\put(1.8,7.85){$\bullet$}
			\put(1.8,9.15){$\bullet$}
			\put(1.5,0.0){$ \Der(L_{5}) $}
			\end{picture}
			
			\item[$ (2) $] The Hasse diagram of $ \Der(L_{2}\times L_{2}) $ is given in \cite[Example 4.21(iii)]{DL}, where $ d_{1}-d_{4} $ are in Example \ref{exa06} (1) and others are the same as $\textnormal {DO}(M_{4}) $ in \cite{DL}. And we can get the Hasse diagram of $ \Der(L_{2}\times L_{3}) $ by Table \ref{t2} in Example \ref{exa06} (2) using python. For details, see the Appendix II listing \ref{a2}.
			\begin{figure}[H]
				\centering
				\subfigure[$ \Der(L_{2}\times L_{2})$]{
					\begin{minipage}[b]{0.45\textwidth}
						\centering
						\includegraphics[width=0.8\textwidth]{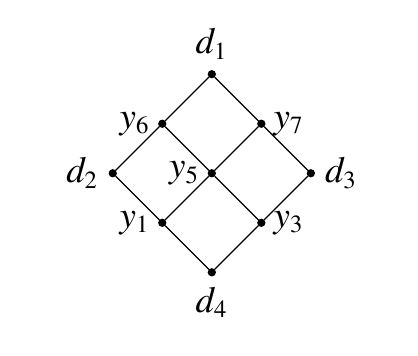} 
					\end{minipage}
				}
				\subfigure[$ \Der(L_{2}\times L_{3}) $]{
					\begin{minipage}[b]{0.45\textwidth}
						\centering
						\includegraphics[width=1.0\textwidth]{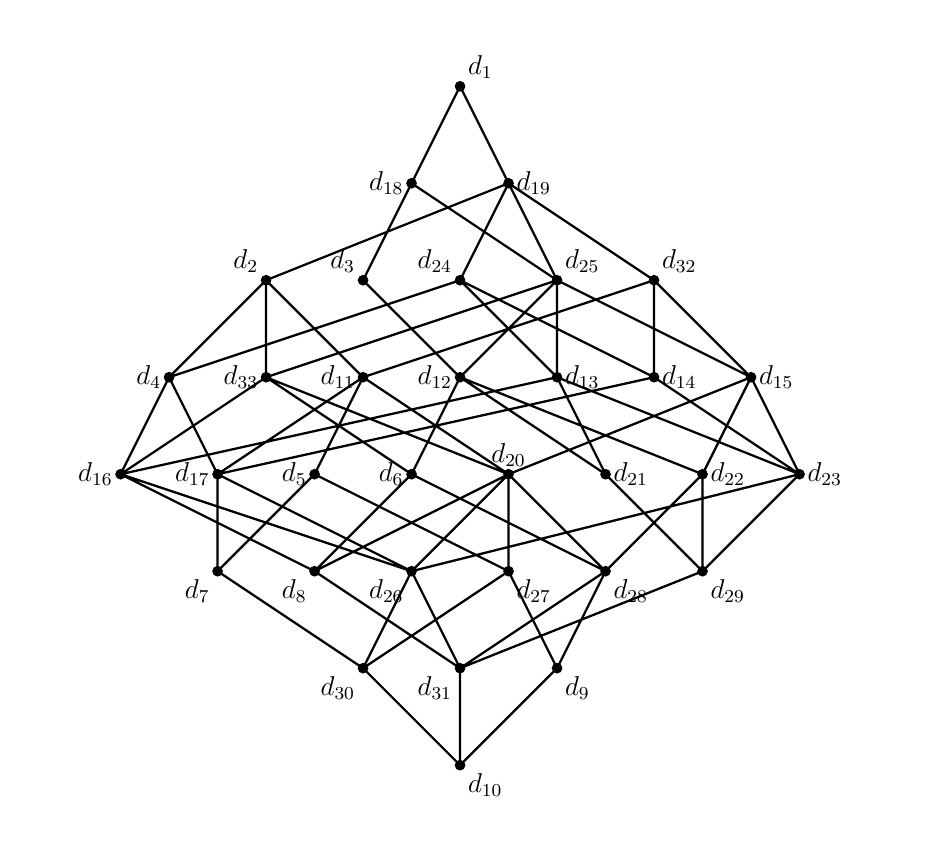}
					\end{minipage}
				}
				\label{02}
			\end{figure}
		\end{enumerate}
	\end{example}
	
	Recall that an MV-algebra $ A $ is \textbf{complete} if its underlying lattice $ \mathbf{L}(A) $ is complete \cite[Definition 6.6.1]{MV4}, that is, for every subset $B$ of $\mathbf{L}(A)$, both $\bigvee B$ and $\bigwedge B$ exist in $\mathbf{L}(A)$. 
	
	Let $\{x_{i}\}_{i\in \Omega}$ be a family elements of $A$ and $x\in A$. If $ \bigvee_{i \in \Omega}x_{i}$ exists, then
	the equality \cite[Chapter V.5]{lattice} holds: 
	\begin{equation}
	x\vee \bigvee_{i \in \Omega}x_{i}=\bigvee_{i \in \Omega}(x\vee x_{i}). 
	\end{equation}

	Let $\{d_{i}\}_{i\in \Omega}$ be a family of operators on a complete MV-algebra $A$. Define  operators  $\bigvee_{i\in \Omega}d_{i}$ and $\bigwedge_{i\in \Omega}d_{i}$ on $A$, respectively, by
	$$(\bigvee_{i\in \Omega}d_{i})(x):=\bigvee_{i\in \Omega}d_{i}(x), \quad
	(\bigwedge_{i\in \Omega}d_{i})(x):=\bigwedge_{i\in \Omega}d_{i}(x)$$
	for any $x\in A$.
	
	
	\begin{lemma}\textnormal{\cite[Lemma 6.6.4]{MV4}}
		Let $A$ be a complete MV-algebra, $x\in A$ and let $\{x_{i}\}_{i\in \Omega}$ be a family elements of $A$.  Then  
		\begin{equation}\label{equ.6}
		x \odot \bigvee_{i \in \Omega} x_{i}=\bigvee_{i \in \Omega}\left(x \odot x_{i}\right).
		\end{equation} 	 
	\end{lemma}
	
	\begin{theorem}\label{thm05}
		Let $A$ be a complete MV-algebra and $\{d_{i}\}_{i\in \Omega}$ be a family elements of $\Der(A)$.  Then the following statements hold:
		\begin{enumerate}
			\item[$ (1) $] $\bigvee_{i\in \Omega}d_{i}\in \Der(A)$.
			\item[$ (2) $]   $(\Der(A), \preceq, \0,\Id)$
			is a complete lattice.
		\end{enumerate}
	\end{theorem}
	\begin{proof}
		
		\vspace{0.15cm} 		
		\noindent$(1)$ For any $x, y\in A$,
		we have
		\begin{eqnarray*}
			(\bigvee_{i\in \Omega}d_{i})(x\odot y)&=& \bigvee_{i\in \Omega}d_{i}(x\odot y)\\
			&\stackrel{(1)}{=}&\bigvee_{i\in \Omega}  ((d_{i}(x)\odot y)\vee (x\odot d_{i}(y)))\\
			&\stackrel{(5)}{=}& (\bigvee_{i\in \Omega}  ((d_{i}(x)\odot y)))\vee \bigvee_{i\in \Omega}(x\odot d_{i}(y))\\
			&\stackrel{(6)}{=}&  ((\bigvee_{i\in \Omega}  d_{i}(x))\odot y))\vee (x\odot \bigvee_{i\in \Omega}d_{i}(y))\\
			&=& ((\bigvee_{i\in \Omega}  d_{i})(x)\odot y)\vee (x\odot (\bigvee_{i\in \Omega}d_{i})(y)),
		\end{eqnarray*}
		and so $\bigvee_{i\in \Omega}d_{i}\in \Der(A)$.
		
		\vspace{0.15cm} 		
		\noindent$(2)$ 
		We shall prove that  $\bigvee_{i\in \Omega}d_{i}$ is the least upper bound of $\{d_{i}\}_{i\in \Omega}$
		in the poset $(\Der(A), \preceq)$.
		Indeed,  firstly, we have 
		$\bigvee_{i\in \Omega}d_{i}\in \Der(A)$ by Item $ (1) $. Secondly, for each $i\in \Omega$, we have
		$d_{i}(x)\leq \bigvee_{i\in \Omega}d_{i}(x)= (\bigvee_{i\in \Omega}d_{i})(x)$ for any $x\in A$ and so
		$d_{i}\preceq \bigvee_{i\in \Omega}d_{i}$. Thus  $\bigvee_{i\in \Omega}d_{i}$ is an upper bound of $\{d_{i}\}_{i\in \Omega}$.
		Finally, let
		$d'\in \Der(A)$ such that  $d_{i}\preceq d'$ for each $i\in \Omega$. Then  $d_{i}(x)\leq d'(x)$  for any $x\in A$, which implies that
		$(\bigvee_{i\in \Omega}d_{i})(x) =\bigvee_{i\in \Omega}d_{i}(x)\leq d'(x)$ and so $ \bigvee_{i\in \Omega}d_{i}\preceq d'$.
		Therefore, we obtain that  $\bigvee_{i\in \Omega}d_{i}$ is
		the least upper bound of $\{d_{i}\}_{i\in \Omega}$
		in the poset $(\Der(A), \preceq)$. Note that $ \bigvee\emptyset=\0 $ and hence
		$(\Der(A), \preceq, \0,\Id)$
		is a complete lattice by \cite[Theorem I.4.2]{bu}.
	\end{proof}
	
	Next we will consider several lattice structure of derivations which are isomorphic to the underlying lattice  $ \mathbf{L}(A) $ of an MV-algebra $A$.
	
	\begin{lemma} \label{lem17}
		Let $A$ be an MV-algebra. Then the following statements hold:
		\begin{enumerate}
			\item[$ (1) $] $d_{u}\vee d_{v}=d_{u\vee v}$ and $d_{u}\wedge d_{v}=d_{u\wedge v}$ for any $u, v\in A$.
			
			\item[$ (2) $] 
			$(\PDer(A), \vee, \wedge, \0,\Id)$ is a sublattice	of $ (\mathrm{O}(A),\preceq) $.
			
			\item[$ (3) $] 	$d\vee d',d\wedge d' \in \IDer(A)$ for any $d, d'\in \IDer(A)$.
			
			\item[$ (4) $] 	$(\IDer(A), \vee, \wedge, \0,\Id)$ is a sublattice	of $ (\mathrm{O}(A),\preceq) $.
		\end{enumerate}
	\end{lemma}
	\begin{proof}
		$ (1) $ Let $u, v\in A$. Then,
		for any $x\in A$, 
		by Lemma \ref{lem02} (6)(7) we have
		$$(d_{u}\vee d_{v})(x)=d_{u}(x)\vee d_{v}(x)=(u\odot x)\vee(v\odot x)=(u\vee v)\odot x =d_{u\vee v}(x),$$
		$$(d_{u}\wedge d_{v})(x)=d_{u}(x)\wedge d_{v}(x)=(u\odot x) \wedge (v\odot x)=(u\wedge v)\odot x=d_{u\wedge v}(x).$$
		Thus $d_{u}\vee d_{v}=d_{u\vee v}$ and $d_{u}\wedge d_{v}=d_{u\wedge v}$.
		
		\vspace{0.15cm} 
		\noindent$ (2) $ It follows immediately from Item (1) that 
		$ \PDer(A) $ is closed under $\vee$ and $\wedge$. So	$(\PDer(A), \vee, \wedge, \0,\Id)$ is a sublattice	of $ (\mathrm{O}(A),\preceq) $, since 
		$\0,\Id\in \PDer(A)$.

		\vspace{0.15cm} 
		\noindent$ (3) $ Let $d, d'\in \IDer(A)$. Then $d(1),d'(1)\in\BA$ and
		$d, d'\in \PDer(A)$ by Proposition \ref{pro04}. Recall that $ \BA $ is a subalgebra of $ A $, since $ d(1),d'(1)\in\BA $, it follows that $ (d\vee d')(1)=d(1)\vee d'(1)= d(1)\oplus d'(1)\in\BA $ by Lemma \ref{lem03} (5). Similarly, $ (d\wedge d')(1)\in\BA $. Moreover, we have $d\vee d', d\wedge d'\in \PDer(A)$ by Item (1). Thus
		$d\vee d', d\wedge d'\in \IDer(A)$.
		
		\vspace{0.15cm} 
		\noindent$ (4) $
		It follows immediately from Item (3) that 
		$ \IDer(A) $ is closed under $\vee$ and $\wedge$. So	$(\IDer(A), \vee, \wedge, \0,\Id)$ is a sublattice	of $ (\mathrm{O}(A),\preceq) $, since 
		$\0,\Id\in \IDer(A)$.	
	\end{proof}
	
	\begin{proposition}\label{pro14}
		Let $A$ be an MV-algebra. Then 
		
		\begin{enumerate}
			\item[$(1)$]
			$(\PDer(A), \vee, \wedge, \0,\Id)$ is a lattice isomorphic to $\mathbf{L}(A)$; and
			
			\item[$(2)$] $(\IDer(A), \vee, \wedge, \0,\Id)$ is a lattice isomorphic to $\BA$.
		\end{enumerate}
		
	\end{proposition}
	\begin{proof}
		$ (1) $ It follows by Lemma \ref{lem17} (2) that 
		$(\PDer(A), \vee, \wedge, \0,\Id)$ is a lattice.
		
		Define a map $g: \PDer(A)\rightarrow \mathbf{L}(A)$ by $g(d_{u})=u$
		for any $d_{u}\in \PDer(A)$. Then $ g $ is a bijection. In fact,  if
		$g(d_{u})=g(d_{v})$, then
		$u=v$, and so $ d_{u}=d_{v} $. Thus $g$ is injective. Also, for each $ u\in A $, there  exists $d_{u}\in\PDer(A)$ such that $ g(d_{u})=u $, so $g$ is surjective. By Lemma \ref{lem17} (1), we have
		$g(d_{u}\vee d_{v})=g(d_{u\vee v})=u\vee v=g(d_{u})\vee g(d_{v})$  and
		$g(d_{u}\wedge d_{v})=g(d_{u\wedge v})=u\wedge v=g(d_{u})\wedge g(d_{v})$.
		Thus $g$ is a lattice isomorphism.
		
		\vspace{0.15cm} 
		\noindent $ (2) $ 
		It follows by Lemma \ref{lem17} (4) that 
		$(\IDer(A), \vee, \wedge, \0,\Id)$ is a lattice.
		
		Define a map $f: \IDer(A)\rightarrow \BA$ by $f(d)=d(1)$
		for any $d\in \IDer(A)$. By Corollary \ref{cor03}, $f$ is a bijection. Also, it is clear that $f(\0)=\0(1)=0$ and
		$f(\Id)=\Id(1)=1$.
		By  Lemma \ref{lem17} (1), we have
		$f(d_{u}\vee d_{v})=f(d_{u\vee v})=u\vee v=f(d_{u})\vee f(d_{v})$  and
		$f(d_{u}\wedge d_{v})=f(d_{u\wedge v})=u\wedge v=f(d_{u})\wedge f(d_{v})$.
		Thus $f$ is a lattice isomorphism.
	\end{proof}
	
	Let $\chi^{(A)}=\{\chi^{(u)}~|~u\in A\}$, where $\chi^{(u)}$ is defined in Corollary \ref{pro09}. We will show that
	$(\chi^{(A)}, \preceq)$ is also a lattice isomorphic to $\mathbf{L}(A)$.

	\begin{lemma}\label{lem18}
		Let $A$ be an MV-algebra and $u, v\in A$. Then the following statements hold:
		\begin{enumerate}
			\item[$ (1) $] $\chi^{(u)}\vee\chi^{(v)}=\chi^{(u\vee v)}$  and
			$\chi^{(u)}\wedge \chi^{(v)}=\chi^{(u\wedge v)}$.
			\item[$ (2) $]  $\chi^{(u)}=\chi^{(v)}$ if and only if $u=v$.
		\end{enumerate}
	\end{lemma}
	\begin{proof}
		$ (1) $ For any $x\in A$, we have
		$$(\chi^{(u)}\vee\chi^{(v)})(x)=\chi^{(u)}(x)\vee\chi^{(v)}(x)=
		\begin{cases}
		u\vee v,  & \textrm{if}~ x=1; \\
		x,  & \textrm{otherwise}
		\end{cases}=\chi^{(u\vee v)}(x)$$
		and
		$$(\chi^{(u)}\wedge\chi^{(v)})(x)=\chi^{(u)}(x)\wedge\chi^{(v)}(x)=
		\begin{cases}
		u\wedge v,  & \textrm{if}~ x=1; \\
		x,  & \textrm{otherwise}
		\end{cases}=\chi^{(u\wedge v)}(x).$$
		Thus $\chi^{(u)}\vee\chi^{(v)}=\chi^{(u\vee v)}$  and
		$\chi^{(u)}\wedge \chi^{(v)}=\chi^{(u\wedge v)}$.
		
		\vspace{0.15cm} 	
		\noindent$(2)$ It is clear that $u=v$ implies $\chi^{(u)}=\chi^{(v)}$. Conversely, if
		$\chi^{(u)}=\chi^{(v)}$, then $u=\chi^{(u)}(1)=\chi^{(v)}(1)=v$.
	\end{proof}
	
	\begin{proposition}\label{pro5.12}
		If $A$ is an MV-algebra, then  $(\chi^{(A)}, \preceq)$ is a sublattice of $(\mathrm{O}(A), \preceq)$ and $(\chi^{(A)}, \preceq)$ is  isomorphic to $\mathbf{L}(A)$.
	\end{proposition}
	\begin{proof}
		Let $u, v\in A$. Then  
		$\chi^{(u)}\vee\chi^{(v)}=\chi^{(u\vee v)}\in \chi^{(A)}$  and
		$\chi^{(u)}\wedge \chi^{(v)}=\chi^{(u\wedge v)}\in \chi^{(A)}$ by Lemma \ref{lem18}. Thus $(\chi^{(A)}, \preceq)$ is a sublattice of $(\mathrm{O}(A), \preceq)$ by Lemma \ref{lem11}.
		
		Define a map $f: \mathbf{L}(A)\rightarrow \chi^{(A)}$ by $f(u)=\chi^{(u)}$
		for any $u\in \mathbf{L}(A)$. By Lemma \ref{lem18}, $f$ is an injective homomorphism. Also, it is clear that $f$ is surjective by the definition of $\chi^{(A)} $. Hence
		$f$ is a lattice isomorphism.
	\end{proof}
	
	Recall that a \textbf{filter} \cite{MV4} of a lattice $L$ is a non-empty subset $F$ of $L$ such that:
	$(i)$  $a, b\in F$ implies $a\wedge b\in F$ and
	$(ii)$ $a\in F$, $c\in L$ and $a\leq c$ imply $c\in F$.
	
	\begin{proposition}
		Let $A$ be an MV-algebra.
		If $(\Der(A), \vee, \wedge, \0,\Id)$
		is a lattice,  then  $\chi^{(A)}$ is a filter of the lattice $\Der(A)$.
	\end{proposition}
	\begin{proof} 
		Assume that $(\Der(A), \vee, \wedge, \0,\Id)$ is a lattice. It is clear that $\chi^{(A)}$ is a non-empty subset of $\Der(A)$ since $ \chi^{(0)}\in\chi^{(A)} $. Also, by Lemma \ref{lem18}, $\chi^{(A)}$ is closed under $\wedge$.
		
		Finally, assume that $d\in \Der(A)$ such that $\chi^{(u)}\preceq d$ for some $u\in A$. Then $A\backslash\{1\}\subseteq \F_{d}(A)$. In fact, for any
		$x\in A\backslash \{1\}$, we have $x=\chi^{(u)}(x)\leq d(x)$ and so $d(x)=x$, since $d(x)\leq x$  by Proposition \ref{pro01} (4).
		It follows that $x\in \F_{d}(A)$ and hence  $A\backslash\{1\}\subseteq \F_{d}(A)$. Consequently, we have
		$d\in \chi^{(A)}$. Therefore, $\chi^{(A)}$ is a filter of the lattice $\Der(A)$.
	\end{proof}

	\section{Discussions}
	~~~~In this paper, we give a detailed algebraic study of $ (\odot,\vee)$-derivations on MV-algebra. There are many different types of derivation on MV, which may lead to more researches and applications. 
	
	We list some questions at the end of this paper. 
	
	1. We have seen in Proposition \ref{pro13} that the relation between the cardinality of MV-algebra $ |A| $ and the cardinality of derivation $ |\Der(A)| $ under small orders. The question is whether we can find the relation when consider larger cardinary $ |\Der(A)| $?
	
	2. In any finite MV-algebra $ A $, we have shown that $ (\Der(A),\preceq, \0,\Id) $ is a lattice in Proposition \ref{pro12} (2). Can we characterize the Hasse diagram of it?
	
	3. In Lemma \ref{lem5.5}, it has been shown that for any MV-chain $ L_{n}(n\geq 2) $, $(\Der(L_{n}),\preceq)$ is a lattice. Naturally, we will ask: for any MV-algebra $A$, is the poset $(\Der(A),\preceq,\0,\Id)$ a lattice?
	
	4. For any two MV-algebras $A$ and $A'$, if $(\Der(A),\preceq,\0,\Id)$ and $(\Der(A'),\preceq,\mathbf{0}_{A'},\textnormal{Id}_{A'})$ are isomorphic lattices, then are $ A $ and $ A' $ isomorphic?
	
	\section*{Declaration}
	~~~~This article does not use any particular data, or human participant. Indeed, the results obtained have been established from the articles cited in the references. However, we remain ready to transmit any information useful for a good understanding of our article.
	
	\textbf{(1) Ethical approval}: We declare that we have complied with the ethical standards for publishing articles in this journal.
	
	\textbf{(2) Funding details}: The work is partially supported by CNNSF (Grants: 12171022, 62250001).
	
	\textbf{(3) Conflict of interest}: The authors have no conflicts of interest to declare that are relevant to the content of this article.
	
	\textbf{(4) Informed Consent}: Not applicable.
	
	\textbf{(5) Authorship contributions}: All authors contributed to this article.

	\section*{Appendix I. Calculation program by Python in Example \ref{exa06} (2)}\label{a01}
	\begin{spacing}{1}
		\begin{lstlisting}[language={Python},caption = {\bf{$\Der(L_{2}\times L_{3})$}.py},label = {a1}]
		ss = ['0', 'a', 'b', 'c', 'd', '1']
		alphabet = {}
		for i in range(len(ss)):
		alphabet[ss[i]] = i
		
		
		def cheng(a, b):
		chenglist = [['0', '0', '0', '0', '0','0'], ['0', '0', 'a', '0', '0', 'a'], ['0', 'a', 'b', '0', 'a','b'], ['0', '0', '0', 'c', 'c', 'c'], ['0', '0', 'a', 'c', 'c', 'd'], ['0', 'a', 'b', 'c', 'd', '1']]
		return chenglist[alphabet[a]][alphabet[b]]
		
		
		def join(a, b):
		joinlist = [['0', 'a', 'b', 'c', 'd', '1'], ['a', 'a', 'b', 'd', 'd', '1'], ['b', 'b', 'b', '1', '1', '1'], ['c', 'd', '1', 'c', 'd', '1'], ['d', 'd', '1', 'd', 'd', '1'], ['1', '1', '1', '1', '1', '1']]
		return joinlist[alphabet[a]][alphabet[b]]
		
		
		sss = []
		for i in ss:
		for j in ss:
		sss.append([i, j])
		
		for a in ['0', 'a']:
		for b in ['0', 'a', 'b']:
		for c in ['0', 'c']:
		for d in ['0', 'a', 'c', 'd']:
		for I in ['0', 'a', 'b', 'c', 'd', '1']:
		mapping = {
		'0': '0', 'a': a, 'b': b, 'c': c, 'd': d, '1': I
		}
		flag = 1
		for i in sss:
		if flag == 1:
		if mapping[cheng(i[0], i[1])] != join(
		cheng(mapping[i[0]], i[1]),
		cheng(i[0], mapping[i[1]])):
		flag = 0
		if flag == 1:
		print(a + b + c + d + I)
		\end{lstlisting}
	\end{spacing}
	\section*{Appendix II. Calculation program by Python in Example \ref{exa07} (2)}\label{a2}
	\begin{spacing}{1}
		\begin{lstlisting}[language={Python},caption = {\bf{Hasse diagram of $\Der(L_{2}\times L_{3})$}.py},label = {a2}]	
		a = [
		'00cc0', '00ccc', '0a000', '0a00a', '0acc0', '0acca', '0accc', '0accd', 'a00a0', 'a0cd0', 'a0cdc', 'aa0a0', 'aa0aa', 'aacd0', 'aacda', 'aacdc', 'aacdd', 'ab000', 'ab00a', 'ab0a0', 'ab0aa', 'ab0ab', 'abcc0', 'abcca', 'abccc', 'abccd', 'abcd0', 'abcda', 'abcdb', 'abcdc', 'abcdd'
		]
		n = len(a[0])
		b = list(i for i in a)
		
		R = ['0a', '0c', '0b', '0d', 'ab', 'ad', 'cd']
		
		
		def leq(a, b):
		if a == b:
		return 0
		for i in range(n):
		if str(a[i]) != str(b[i]) and str(a[i]) + str(b[i]) not in R:
		return 0
		return 1
		
		
		def maximal(a):
		max = []
		for i in a:
		b = set(leq(i, j) for j in a)
		if 1 not in b:
		max.append(i)
		return max
		
		
		def minimal(a):
		min = []
		for i in a:
		b = set(leq(j, i) for j in a)
		if 1 not in b:
		min.append(i)
		return min
		
		
		while b != []:
		print(minimal(b))
		b = [i for i in b if not i in minimal(b)]
		print()
		while a != []:
		print(maximal(a))
		a = [i for i in a if not i in maximal(a)]
		
		\end{lstlisting}
	\end{spacing}

\end{document}